\documentclass[12pt]{elsarticle}
\usepackage{MnSymbol,longtable,amsmath,amsthm,bbold}
\usepackage[all]{xy}
\usepackage[active]{srcltx}
\usepackage[T1]{fontenc}
\DeclareMathOperator{\codim}{codim}
\DeclareMathOperator{\diag}{diag}
\DeclareMathOperator{\rank}{rank}
\DeclareMathOperator{\gl}{GL}
\DeclareMathOperator{\SL}{SL}
\DeclareMathOperator{\st}{ig}
\DeclareMathOperator{\Or}{O}
\DeclareMathOperator{\Sp}{Sp}

\newcommand{\C}{\mathbb C}

\newcommand{\ca}{_{\text{\rm can}}}

\newcommand{\mat}[1]{%
\left[\begin{smallmatrix} #1
 \end{smallmatrix}\right]}

\newcommand{\set}[2]{%
\left\{\left.#1\right|#2\right\}}

\renewcommand{\le}{\leqslant}

\newtheorem{theorem}{Theorem}[section]
\newtheorem{proposition}[theorem]{Proposition}
\newtheorem{lemma}[theorem]{Lemma}
\newtheorem{corollary}[theorem]{Corollary}

\theoremstyle{definition}
\newtheorem{definition}[theorem]{Definition}

\theoremstyle{remark}
\newtheorem{remark}[theorem]{Remark}

\begin{document}

\title{Change of
the congruence canonical form of 2-by-2
and 3-by-3 matrices under perturbations
and bundles of matrices under
congruence}

\author[dmi]{Andrii Dmytryshyn}
\ead{andrii@cs.umu.se}
\address[dmi]{Department of Computing Science
and HPC2N, Ume{\aa} University, Sweden}

\author[fut]{Vyacheslav Futorny}
\ead{futorny@ime.usp.br}

\address[fut]{Department of Mathematics, University
of S\~ao Paulo, Brazil}

\author[dmi]{Bo K\r{a}gstr\"{o}m}
\ead{bokg@cs.umu.se}

\author[kli]{Lena~Klimenko}
\ead{e.n.klimenko@gmail.com}
\address[kli]{National Technical
University of Ukraine ``Kyiv
Polytechnic Institute'', Kiev, Ukraine}

\author[ser]{Vladimir V.
Sergeichuk\corref{cor}}
\ead{sergeich@imath.kiev.ua}
\address[ser]{Institute of Mathematics,
Kiev, Ukraine}
\cortext[cor]{Corresponding author; published in Linear Algebra Appl. 469 (2015) 305-334}

\begin{abstract}
We construct the Hasse diagrams $G_2$
and $G_3$ for the closure ordering on
the sets of congruence classes of
$2\times 2$ and $3\times 3$ complex
matrices. In other words, we construct
two directed graphs whose vertices are
$2\times 2$ or, respectively, $3\times
3$ canonical matrices under congruence
and there is a directed path from $A$
to $B$ if and only if $A$ can be
transformed by an arbitrarily small
perturbation to a matrix that is
congruent to $B$.

A bundle of matrices under congruence
is defined as a set of square matrices
$A$ for which the pencils $A+\lambda
A^T$ belong to the same bundle under
strict equivalence. In support of this
definition, we show that all matrices
in a congruence bundle of $2\times 2$
or $3\times 3$ matrices have the same
properties with respect to
perturbations. We construct the Hasse
diagrams $G_2^{\rm B}$ and $G_3^{\rm
B}$ for the closure ordering on the
sets of congruence bundles of $2\times
2$ and, respectively, $3\times 3$
matrices. We find the isometry groups
of $2\times 2$ and $3\times 3$
congruence canonical matrices.
\end{abstract}

 \maketitle
\section{Introduction}
\label{introd}

We study how small perturbations of a
$2\times 2$ or $3\times 3$ complex
matrix can change its congruence
canonical form.

Two complex matrices $A$ and $B$ are
said to be \emph{congruent}  if
$S^TAS=B$ for a nonsingular $S$. This
is an equivalence relation; its
equivalence classes are called
\emph{congruence classes}. In Section
\ref{s2} we construct the closure
graphs $G_2$ and $G_3$, which are
defined for any natural $n$ as follows.

\begin{definition}\label{kur}
The \emph{closure graph  $G_n$} for
congruence classes of $n\times n$
complex matrices is the directed graph
in which
 each vertex $v$ represents
      in a one-to-one manner a
      congruence class $C_v$ of
      $n\times n$ matrices, and
  there is a directed path
      from a vertex $v$ to a vertex
      $w$ if and only if one (and hence each) matrix
      from $C_v$ can be transformed
     to a matrix from
      $C_w$  by an arbitrarily small
      perturbation.
\end{definition}

The graph $G_n$ is the Hasse diagram of
the set of congruence classes of
$n\times n$ matrices with the following
partial order: $a\preccurlyeq b$ if $a$
is contained in the closure of $b$.
Thus, the graph $G_n$ shows how the
congruence classes relate to each other
in the affine space of $n\times n$
matrices.

Since each $n\times n$ matrix is
uniquely represented in the form $P+Q$
in which $P$ is symmetric and $Q$ is
skew-symmetric, $G_n$ is also the
closure graph for congruence classes of
$n\times n$ symmetric/skew-symmetric
matrix pencils $P+\lambda Q$.

Each congruence class contains exactly
one canonical matrix for congruence,
and so it is convenient to represent
the congruence classes by their
canonical matrices. We use the
congruence canonical matrices
$A_{\text{can}}$ constructed by Horn
and Sergeichuk \cite{hor-ser_transp}.
We also use the \emph{miniversal
deformation} of $A_{\text{can}}$ given
by Dmytryshyn, Futorny, and Sergeichuk
\cite{f-ser}; that is, a simple normal
form to which all matrices close to
$A_{\text{can}}$ can be reduced by
congruence transformations that
smoothly depend on their entries.

The closure graph for *congruence
classes of $2\times 2$ complex matrices
was constructed by Futorny, Klimenko,
and Sergeichuk \cite{f_gr*cong}.

Unlike perturbations of matrices under
congruence and *congruence,
perturbations of matrices under
similarity and of matrix pencils have
been much studied. For a given matrix
$A$, Boer and Thijsse \cite{den-thi}
and, independently, Markus and Parilis
\cite{mar-par} described the set of all
Jordan canonical matrices $J$ such that
in each neighborhood of $A$ there
exists a matrix whose Jordan canonical
form is $J$. Their description was
extended to Kronecker's canonical forms
of pencils by Pokrzywa \cite{pok}.

Arnold \cite[\S\,5.3]{arn} defines a
\emph{bundle of matrices under
similarity} as a set of all matrices
having the same Jordan type: matrices
$A$ and $B$ have the same \emph{Jordan
type} if there is a bijection from the
set of distinct eigenvalues of $A$ to
the set of distinct eigenvalues of $B$
that transforms the Jordan canonical
form of $A$ to the Jordan canonical
form of $B$. For example, the matrices
\begin{equation}\label{lic}
J_3(0)\oplus J_2(0)\oplus
J_5(1)\quad\text{and}\quad J_3(2)\oplus J_2(2)\oplus
J_5(-3)
\end{equation}
belong to the same bundle (we denote by
$J_n(\lambda )$ the $n\times n$ upper
triangular Jordan block with eigenvalue
$\lambda$). All matrices of a bundle
have similar properties; for example,
its Jordan matrices $J$ have the same
set  $\{X\,|\, JX=XJ\}$ of commuting
matrices.

Note that the closure graph for bundles
of $n\times n$ matrices under
similarity has a finite number of
vertices; moreover, it is in some sense
more informative than the closure graph
for similarity classes. For example,
one cannot see from the latter graph
that each neighborhood of $J_n(\lambda
)$ contains a matrix with $n$ distinct
eigenvalues (since there is no diagonal
matrix whose similarity class has a
nonzero intersection with each
neighborhood of $J_n(\lambda )$). But
the closure graph for bundles has a
directed path from the bundle
containing $J_n(\lambda )$ to the
bundle of all matrices with $n$
distinct eigenvalues.

The bundles of matrix pencils are
defined in the same way via Kronecker's
canonical form of pencils (see Section
\ref{ghy}). Edelman, Elmroth, and
K\r{a}gstr\"{o}m \cite{kag2} developed
a comprehensive theory of closure
relations for similarity classes of
matrices, for strict equivalence
classes of matrix pencils, and for
their bundles. The software StratiGraph
\cite{e-j-k} constructs their closure
graphs.  The closure graph for $2\times
3$ matrix pencils was constructed and
studied by Elmroth and K\r{a}gstr\"{o}m
\cite{El-ka}.

The definition of bundles of matrices
under congruence is not so evident. One
could define these bundles via the
congruence canonical form (see
Definition \ref{mos}) by analogy with
bundles of matrices under similarity
and bundles of matrix pencils. But
unlike the Jordan and Kronecker
canonical forms, the perturbation
behavior of a congruence canonical
matrix with parameters depends on the
values of its parameters (see Remark
\ref{lrp}). Moreover, one can obtain
another partition into bundles using
another congruence canonical form (say,
the tridiagonal canonical form
\cite[Theorem 1.1]{f-h-s_trid}). We
define congruence bundles as follows.

\begin{definition}\label{kuq1}
Two square matrices $A$ and $B$ belong
to the same \emph{congruence bundle} if
and only if
 the pencils $A+\lambda A^T$
and $B+\lambda B^T$ belong to the same
strict equivalence bundle.
\end{definition}

This definition is based on the
remarkable fact: two $n\times n$
matrices $A$ and $B$ are congruent if
and only if the pencils $A+\lambda A^T$
and $B+\lambda B^T$ are strictly
equivalent (see Lemma \ref{leq}).

An informal introduction to
perturbations of matrices determined up
to similarity, congruence, or
*congruence is given by Klimenko and
Sergeichuk \cite{k-s_deform}.
Miniversal deformations of matrices
under similarity, congruence, and
*congruence and of matrix pencils under
strict equivalence are given in
\cite{arn,arn1,f-ser,def-sesq,kag,gar_ser,%
k-s_triang}.

The term ``congruence orbit'' is often
used instead of ``congruence class''
\cite{ter-dor,dm_rep}. The problem that
we consider can be called ``the
stratification of orbits and bundles of
matrices under congruence'' by analogy
with the stratification of orbits and
bundles of matrices under similarity
and of matrix pencils
\cite{kag2,e-j-k2,joh1,joh}. The
stratification theory for
skew-symmetric matrix pencils has
recently been developed in
\cite{ssstrat, sscodim}.

The paper is organized as follows. In
Section \ref{s2} we construct the
closure graphs $G_2$ and $G_3$ for
congruence classes of $2\times 2$ and
$3\times 3$ matrices. In Section
\ref{s3} we construct the closure
graphs $G^{\rm B}_2$ and $G^{\rm B}_3$
for congruence bundles of $2\times 2$
and $3\times 3$ matrices. In Sections
\ref{ssk} and \ref{sss} we prove the
main theorems.

In Section \ref{kjr} we give arguments
in favor of Definition \ref{kuq1}. We
show that each congruence bundle of
$2\times 2$ or $3\times 3$ matrices can
be described as a set of all matrices
with similar properties with respect to
perturbations and with the same number
of indecomposable direct summands
(which is not true for bundles defined
via canonical matrices for congruence;
see Definition \ref{mos}). Thus, one
obtains the same closure graphs $G^{\rm
B}_2$ and $G^{\rm B}_3$ using any other
\emph{good} definition of congruence
bundles. In Section \ref{kjtf} we give
the list of isometry groups of bilinear
spaces $(\mathbb C^n,\mathcal A)$ in
which $\cal A$ is a bilinear form on
$\mathbb C^n$ given by a congruence
canonical matrix and $n=2$ or $3$.

All matrices that we consider are over
the field of complex numbers.

\section{Closure graphs for congruence
classes}\label{s2}

We use the following canonical form
under congruence (in which $I_m$ is the
$m \times m$ identity matrix).

\begin{proposition}
[{\cite[Theorem 4.5.25]{hor-jo}}]
\label{le1} Every square
      complex matrix is congruent
      to a direct sum, determined
      uniquely up to permutation of
      summands, of matrices of the
      form
\begin{equation}
\label{can}
H_{2m}(\lambda):=
\begin{bmatrix}0&I_m\\
J_m(\lambda) &0
\end{bmatrix},
 \qquad
\Gamma_n:= \begin{bmatrix} 0&&&&
\udots
\\&&&-1&\udots
\\&&1&1\\ &-1&-1& &\\
1&1&&&0
\end{bmatrix},
 \qquad
J_k(0),
  \end{equation}
in which $\lambda$ is a nonzero complex
number that is determined up to
replacement by $\lambda^{-1}$ and
satisfies $\lambda\ne (-1)^{m+1}$, and
$\Gamma _n$ is $n$-by-$n$.
\end{proposition}

This canonical form and an analogous
canonical form for *congruence were
obtained in \cite{hor-ser_transp} based
on \cite[Theorem 3]{ser_izv} (see also
\cite{hor-ser_anyf}); a direct proof
that the forms are canonical is given
in \cite{hor-ser_regul, hor-ser_can}.

For each $A\in{\mathbb C}^{n\times n}$
and a matrix $X\in{\mathbb C}^{n\times
n}$ with small entries,
\[
(I+X)^TA(I+X)
=A+\underbrace{X^T A+ AX}
_{\text{small}}
+\underbrace{X^T AX}
_{\text{very small}}
\]
and so the congruence class of $A$ in a
small neighborhood of $A$ can be
obtained by a very small deformation of
the affine matrix space \[\{A+ X^T A+
AX\,|\,X\in{\mathbb C}^{n\times n}\}.\]
(By the local Lipschitz property
\cite{rodm}, if $A$ and $B$ are close
to each other and $B=S^TAS$ with a
nonsingular $S$, then $S$ can be taken
near $I_n$.) The vector space
\begin{equation}\label{msiy}
T(A):=\{X^TA+AX\,|\,X\in{\mathbb
C}^{n\times n}\}
\end{equation}
is the tangent space to the congruence
class of $A$ at the point $A$. The
numbers
\begin{equation}\label{kow}
\dim_{\mathbb C} T(A),\qquad
\codim_{\mathbb C} T(A):=
n^2-\dim_{\mathbb C} T(A)
\end{equation}
are called the \emph{dimension} and,
respectively, \emph{codimension} of the
congruence class of $A$.

\begin{theorem}[proved in Section
\ref{s_pr1}]\label{the1} The closure
graph $G_2$ for congruence classes of\/
$2\times 2$ matrices is given in Figure
\ref{kib5}.
 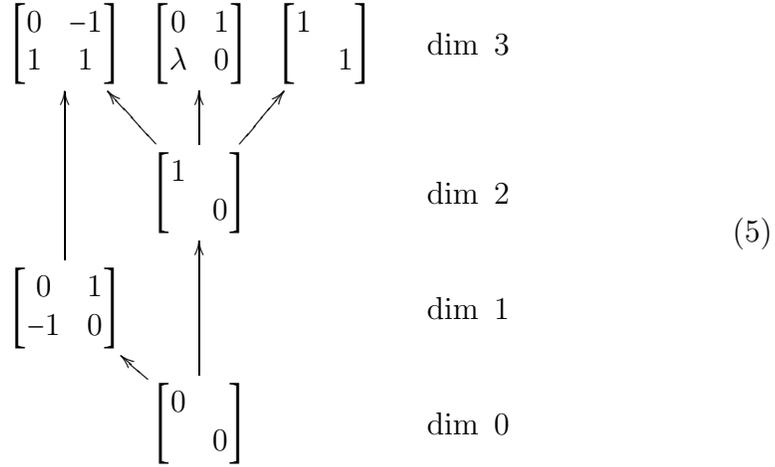
\begin{figure}[hbt]
\begin{equation}\label{g1}
\begin{split}
\xymatrix@R=7pt@C=8pt{%
{\begin{bmatrix} 0&-1\\ 1&1
 \end{bmatrix}}
    &
{\begin{bmatrix} 0&1\\ \lambda &0
\end{bmatrix}}
   &{\begin{bmatrix} 1&\\
&1 \end{bmatrix}}
 &
  *[r]{\dim\ 3}
    \\ \\
  &{\begin{bmatrix}
1&\\ &0
 \end{bmatrix}}\ar[uu] \ar[luu]
\ar[uur]&
 &
 *[r]{\dim\ 2}
    \\
{\begin{bmatrix} 0&1\\
-1&0
 \end{bmatrix}}\ar[uuu]& & &
*[r]{\dim\ 1}
      \\
&{\begin{bmatrix} 0&\\
&0
 \end{bmatrix}}
\ar[ul]\ar[uu] & &
*[r]{\dim\ 0}
} \end{split}
\end{equation}
\caption{\small The closure
graph $G_2$ for congruence classes of\/
$2\times 2$ matrices, in which
$\lambda\ne\pm 1$, and each nonzero $\lambda
$ is determined up to
replacement by $\lambda^{-1} $.} \label{kib5}
\end{figure}
Each congruence class is represented by
its canonical matrix, which is a direct
sum of blocks of the form \eqref{can}
$($the zero entries outside of these
blocks are not shown$)$. The graph is
infinite:
$\left[\begin{smallmatrix} 0&1\\
\lambda &0
\end{smallmatrix}\right]$ represents the infinite
set of vertices indexed by $\lambda
\in\mathbb C \smallsetminus\{-1,1\}$
$($provided that each nonzero $\lambda$
is determined up to replacement by
$\lambda ^{-1})$ with arrows
$\left[\begin{smallmatrix} 1&\\ &0
\end{smallmatrix}\right]\to
\left[\begin{smallmatrix} 0&1\\
\lambda &0
\end{smallmatrix}\right]$ to each of these vertices.
The congruence classes of canonical
matrices located at the same horizontal
level in $G_2$ have the same dimension,
which is indicated to the right.
\end{theorem}

For example, the graph $G_2$ in
\eqref{g1} shows that an arbitrarily
small
      neighborhood of
$\left[\begin{smallmatrix} 0&1\\
-1 &0
\end{smallmatrix}\right]$ contains
matrices with congruence canonical
forms
$\left[\begin{smallmatrix} 0&1\\
-1 &0
\end{smallmatrix}\right]$ and
$\left[\begin{smallmatrix} 0&-1\\
1 &1
\end{smallmatrix}\right],$
but for any other canonical matrix
$A_{\text{can}}$ there is a
neighborhood of $\left[\begin{smallmatrix} 0&1\\
-1 &0
\end{smallmatrix}\right]$ without
matrices with canonical form
$A_{\text{can}}$.

\begin{theorem}[proved in Section
\ref{s_pr2}]\label{the2} The closure
graph $G_3$ for congruence classes of\/
$3\times 3$ matrices is given in Figure
\ref{fig2}.
\begin{figure}[hbt]
\begin{equation}\label{g3}
\begin{split}
\xymatrix@C=10pt@R=14pt{
{\left[\begin{smallmatrix} 0&-1&\\
1&1&\\&&1
\end{smallmatrix}\right]}&&
{\left[\begin{smallmatrix}
0&1&\\ {\mu} &0&\\&&1
\end{smallmatrix}\right]} &&
{\left[\begin{smallmatrix}
0&0&1\\ 0&-1&-1\\1&1&0
\end{smallmatrix}\right]}&&&\dim\, 8
    \\
&& {\left[\begin{smallmatrix}
0&1&0\\ 0&0&1\\0&0&0
\end{smallmatrix}\right]} \ar[urr]\ar[ull] \ar[u]&&
              &&&\dim\, 7
\\
{\left[\begin{smallmatrix}
0&1&\\ -1&0&\\&&1
\end{smallmatrix}\right]}\ar[uu]&&&&
{\left[\begin{smallmatrix}
1&&\\ &1&\\&&1
\end{smallmatrix}\right]}\ar[uu]
                           &&&\dim\, 6
\\
{\left[\begin{smallmatrix}
0&-1&\\ 1&1&\\&&0
\end{smallmatrix}\right]}
\ar[uurr]\ar[u]
&
&{\left[\begin{smallmatrix}
0&1&\\ {\lambda}
&0&\\&&0
\end{smallmatrix}\right]}\ar[uu]
&&{\left[\begin{smallmatrix}
1&&\\ &1&\\&&0
\end{smallmatrix}\right]} \ar[uull]\ar[u]
                        &&&\dim\, 5
         \\
{\left[\begin{smallmatrix}
0&1&\\ -1&0&\\&&0
\end{smallmatrix}\right]}
\ar[u] &&
{\left[\begin{smallmatrix}
1&&\\ &0&\\&&0
\end{smallmatrix}\right]}\ar[u]
\ar[ull]
\ar[urr]&&
                         &&&\dim\, 3
            \\
&&{\left[\begin{smallmatrix}
0&&\\ &0&\\&&0
 \end{smallmatrix}\right]}\ar[ull]\ar[u]&&
                 &&&\dim\, 0
      }
\end{split}
\end{equation}
\caption{\small The closure graph $G_3$
for congruence classes of\/ $3\times 3$
matrices, in which
$\lambda,\mu\ne\pm 1$, and nonzero $\lambda
$ and $\mu $ are determined up to
replacements by $\lambda^{-1} $ and
$\mu^{-1}$.} \label{fig2}
\end{figure}
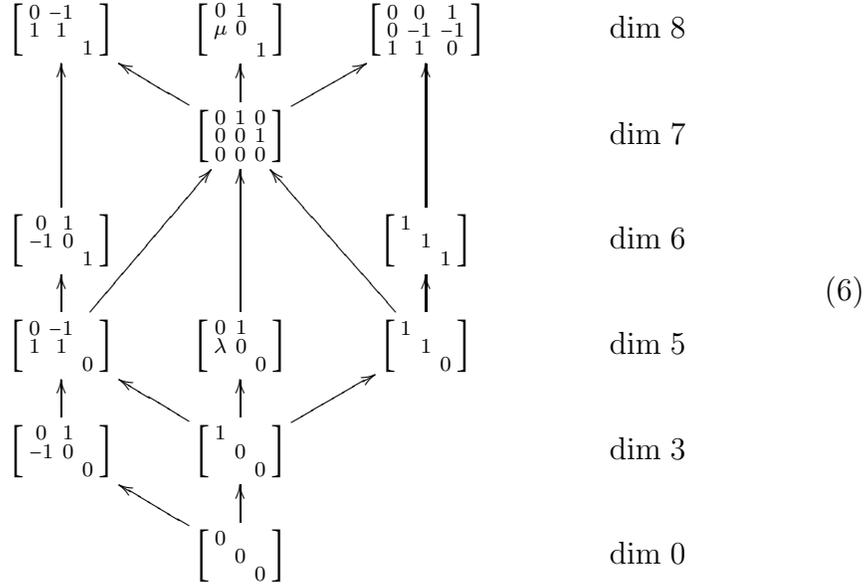
The congruence classes that correspond
to the vertices are represented by
their $3\times 3$ canonical matrices
under congruence.
 The graph is infinite:
\[\left[\begin{smallmatrix} 0&1&\\
\lambda &0&\\&&0
 \end{smallmatrix}\right]\quad\text{and}
\quad\left[\begin{smallmatrix} 0&1&\\ \mu
&0&\\&&1
 \end{smallmatrix}\right]\]
represent the infinite sets of vertices
indexed by $\lambda, \mu \in\mathbb C
\smallsetminus\{-1,1\}$ provided that
each nonzero $\lambda$ is determined up
to replacement by $\lambda ^{-1}$ and
each nonzero $\mu$ is determined up to
replacement by $\mu ^{-1}$. The
congruence classes with vertices on the
same horizontal level have the same
dimension, which is indicated to the
right.
\end{theorem}

\begin{remark}\rm
Let $M$ be a $2\times 2$ or $3\times 3$
canonical matrix for congruence.

\begin{itemize}
\item Let $N$ be another canonical
    matrix for congruence of the
    same size. Each neighborhood of
$M$ contains a matrix whose
    congruence canonical form is
$N$ if and only if there is a
directed path from $M$ to $N$ in
 $G_2$ or $G_3$  (if $M=N$ then
    there is always the ``lazy''
    path of length $0$ from $M$ to
    $M$).

  \item The closure of the
      congruence class of $M$ is
      equal to the union of the
      congruence classes of all
      canonical matrices $N$ such
      that there is a directed path
      from $N$ to $M$ (if $N=M$
      then there is always the
      ``lazy'' path).
\end{itemize}
\end{remark}

\section{Closure graphs for congruence
bundles} \label{s3}

The following lemma, which is another
version of Proposition \ref{le1}, is
used for constructing closure graphs
for congruence.

\begin{lemma}
\label{thejr} Every square complex
matrix is
      congruent to a direct sum,
      determined uniquely up to
      permutation of summands, of
      matrices of the form
\begin{equation}
\label{cjy}
H_{2m}(\lambda),
 \qquad
\Gamma _n,
 \qquad
J_{2r-1}(0),
  \end{equation}
in which
$(-1)^{m+1}\ne\lambda\in\mathbb C$,
each nonzero $\lambda$ is determined up
to replacement by $\lambda^{-1}$, and
$m,n,r=1,2,\dots$\,.
\end{lemma}

\begin{proof}
The lemma follows from Proposition
\ref{le1} since by \cite[p.
213]{hor-ser_anyf} or \cite[p.
492]{ser_izv} $J_{k}(0)$ is
permutationally congruent to
\begin{equation}\label{kuf}
\begin{bmatrix}
0&R_r^T\\L_r&0
\end{bmatrix}\text{ if }k =2r-1,\qquad
H_{2r}(0)=\begin{bmatrix}
0&I_r\\J_r(0)&0
\end{bmatrix}\text{ if }k=2r,
\end{equation}
in which
\begin{equation}\label{1.4}
L_r:=\begin{bmatrix}
1&0&&0\\&\ddots&\ddots&\\0&&1&0
\end{bmatrix},\quad
R_r:=\begin{bmatrix}
0&1&&0\\&\ddots&\ddots&\\0&&0&1
\end{bmatrix}
\end{equation}
are $(r-1)\times r$ matrices (thus,
$L_1=R_1$ is the $0\times 1$ matrix).
\end{proof}

The congruence bundles (see Definition
\ref{kuq1}) are described in the
following theorem. By Lemma
\ref{thejr}, each square matrix $A$ is
congruent to a direct sum
\begin{equation}\label{lmio}
{\cal G}\oplus {\cal H}_1(\lambda_1)
\oplus\dots\oplus {\cal H}_t(\lambda_t),
\quad \lambda_i\ne \lambda _j,1/\lambda _j
\text{ if }
i\ne j, \text{ all }\lambda_i
\in\mathbb C\setminus\{1,-1\},
\end{equation}
in which ${\cal G}$ is a direct sum of
matrices of the form $H_{2m}((-1)^m)$,
$\Gamma _n$, $J_{2r-1}(0)$, and each
${\cal H}_i(\lambda_i)$ is a direct sum
of matrices of the form
$H_{2m}(\lambda_i)$ with $\lambda_i\ne
\pm 1$. (The condition $\lambda_i\ne
1/\lambda _j$ in \eqref{lmio} is
obviously satisfied if $\lambda _j=0$.)

\begin{theorem}[proved in Section  \ref{ghy}]
\label{theorm} Let a square matrix $A$
be congruent to \eqref{lmio}. Then the
congruence bundle of $A$ consists of
all matrices that are congruent to
matrices
\begin{equation}\label{kdr}
{\cal G}\oplus {\cal H}_1(\mu_1)
\oplus\dots\oplus {\cal H}_t(\mu_t),
\quad \mu_i\ne \mu _j,1/\mu _j\text{ if }
i\ne j,\text{ all }\mu_i
\in\mathbb C\setminus\{1,-1\}
\end{equation}
with the same ${\cal G},{\cal H}_1,
\dots, {\cal H}_t.$
\end{theorem}

 \begin{corollary}\label{jjy}
The congruence
      bundles that contain the
      canonical blocks \eqref{can}
      are the following:
\begin{itemize}
  \item the congruence classes of
      $H_{2m}((-1)^m)$, $\Gamma
      _n$, and $J_{2r-1}(0)$,

  \item  the union of the
      congruence class of
      $J_{2m}(0)$ and the
      congruence classes of all
      $H_{2m}(\lambda)$ with
      $\lambda\ne 0,\pm 1$ and
$\lambda$ is determined up to
replacement by $\lambda^{-1}$
\end{itemize}
for each $m,n,r=1,2,\dots$\,.
\end{corollary}

\begin{remark}\label{lrp}
The matrices $H_{2m}((-1)^m)$ and
$H_{2m}(\lambda)$ with $\lambda\ne
0,\pm 1$ properly belong to distinct
bundles because $ \left[
\begin{smallmatrix} 0&1\\
-1 &0
 \end{smallmatrix}\right]$ and $\left[
\begin{smallmatrix} 0&1\\
\lambda &0
 \end{smallmatrix}\right]$ with $\lambda
 \ne \pm 1$ have distinct properties
 with respect to perturbations, which
 is illustrated by the closure graph $G_2$ in \eqref{g1}.
\end{remark}

Similar to the closure graph $G_n$ for
congruence classes, we denote by
$G_n^{\rm{B}}$ the closure graph of
congruence bundles.

\begin{theorem}[proved in Section  \ref{gep}]
\label{the1_b}
\begin{itemize}
  \item[{\rm (a)}] The closure
      graph $G_2^{\rm{B}}$ for
      congruence bundles of\/
      $2\times 2$ matrices is given
      in Figure
\ref{fig3}.
\begin{figure}[hbt]
\begin{equation}\label{g2}
\begin{split}
\xymatrix@R=1pt@C=20pt{
{\phantom{{}_{\lambda \ne\pm 1}}\left\{\begin{bmatrix}
0&1\\ \lambda &0
 \end{bmatrix}\right\}_{\lambda \ne\pm 1}}
 &&&\dim\:4
\\    \\ \\
{\begin{bmatrix} 0&-1\\ 1&1
 \end{bmatrix}}\ar[uuu]
    &
{\begin{bmatrix} 1&\\
&1
\end{bmatrix}}\ar[uuul]
   &&\dim\:3
 \\   \\ \\
  &{\begin{bmatrix}
1&\\ &0
 \end{bmatrix}}\ar[uuu] \ar[luuu]
&&\dim\:2
    \\
{\begin{bmatrix} 0&1\\
-1&0
 \end{bmatrix}}\ar[uuu]& &
 &\dim\:1
      \\
&{\begin{bmatrix} 0&\\
&0
 \end{bmatrix}}
\ar[ul]\ar[uu] &&\dim\:0
} \end{split}
\end{equation}
\caption{\small The closure
      graph $G_2^{\rm{B}}$ for congruence bundles
      of\/ $2\times 2$ matrices.}
      \label{fig3}
\end{figure}
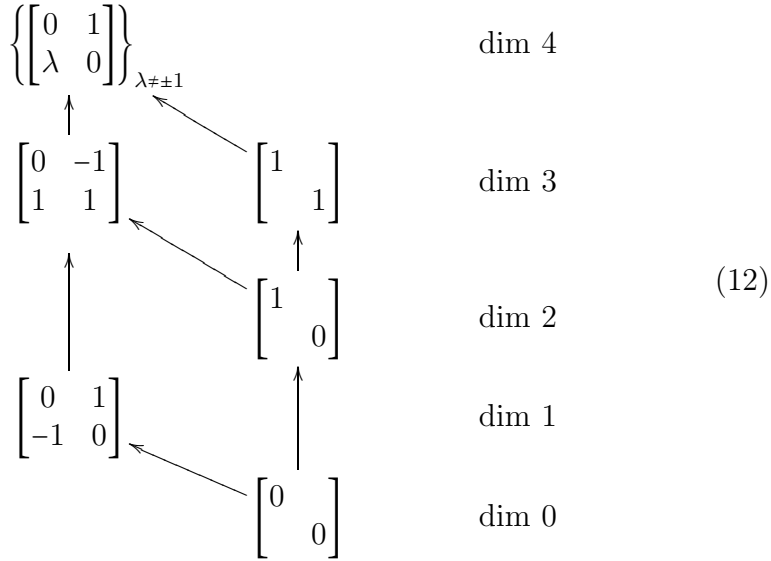
The vertex $\left\{ \left[
\begin{smallmatrix} 0&1\\
\lambda &0
 \end{smallmatrix}\right]
 \right\}_{\lambda \ne \pm 1}$
represents the bundle that consists
of all matrices whose congruence
canonical forms are
$\left[\begin{smallmatrix} 0&1\\
\lambda &0
 \end{smallmatrix}\right]$ with $\lambda
 \ne \pm 1$; each nonzero $\lambda$
is determined up to replacement by
$\lambda^{-1}$. The other vertices
are congruence canonical matrices;
the corresponding bundles coincide
with their congruence classes.

  \item[{\rm (b)}] The closure
      graph $G_3^{\rm{B}}$ for
      congruence bundles of\/
      $3\times 3$ matrices is given
      in Figure
\ref{fig4}.
\begin{figure}[hbt]
\begin{equation}\label{g4}
\begin{split}
\xymatrix@C=0pt@R=14pt{
&&
{\phantom{{}_{\lambda \ne\pm 1}}\left\{\left[\begin{smallmatrix}
0&1&\\ {\mu} &0&\\&&1
\end{smallmatrix}\right]\right\}_{\mu \ne\pm 1}} &&
&\quad&{\dim\, 9}
\\
{\left[\begin{smallmatrix}
0&-1&\\ 1&1&\\&&1
\end{smallmatrix}\right]}\ar[urr]&&&&
{\left[\begin{smallmatrix}
0&0&1\\ 0&-1&-1\\1&1&0
\end{smallmatrix}\right]}\ar[ull]&
\quad&{\dim\, 8}
    \\
&& {\left[\begin{smallmatrix}
0&1&0\\ 0&0&1\\0&0&0
\end{smallmatrix}\right]} \ar[urr]\ar[ull]&& &
\quad&{\dim\, 7}
    \\
{\left[\begin{smallmatrix}
0&1&\\ -1&0&\\&&1
\end{smallmatrix}\right]}\ar[uu]&&
{\phantom{{}_{\lambda \ne\pm 1}}\left\{\left[\begin{smallmatrix}
0&1&\\ {\lambda}
&0&\\&&0
\end{smallmatrix}\right]
\right\}_{\lambda \ne\pm 1}}\ar[u]&&
{\left[\begin{smallmatrix}
1&&\\ &1&\\&&1
\end{smallmatrix}\right]}\ar[uu]&
\quad&{\dim\, 6}\\
{\left[\begin{smallmatrix}
0&-1&\\ 1&1&\\&&0
\end{smallmatrix}\right]}
\ar[u]\ar[urr]&
&
&&{\left[\begin{smallmatrix}
1&&\\ &1&\\&&0
\end{smallmatrix}\right]} \ar[ull]\ar[u] &
\quad&{\dim\, 5}
         \\
{\left[\begin{smallmatrix}
0&1&\\ -1&0&\\&&0
\end{smallmatrix}\right]}
\ar[u] &&
{\left[\begin{smallmatrix}
1&&\\ &0&\\&&0
\end{smallmatrix}\right]}
\ar[ull]
 \ar[urr]&&
 &
\quad&{\dim\, 3}\\
&&{\left[\begin{smallmatrix}
0&&\\ &0&\\&&0
 \end{smallmatrix}\right]}\ar[ull]\ar[u]&&
 &
\quad&{\dim\, 0}
      }
\end{split}
\end{equation}
\caption{\small\it The closure graph $G_3^{\rm{B}}$ for
congruence bundles of
${3\times 3}$ matrices, in which
$\lambda,\mu\ne\pm 1$, and nonzero $\lambda
$ and $\mu $ are determined up to
replacements by $\lambda^{-1} $ and
$\mu^{-1}$.} \label{fig4}
\end{figure}
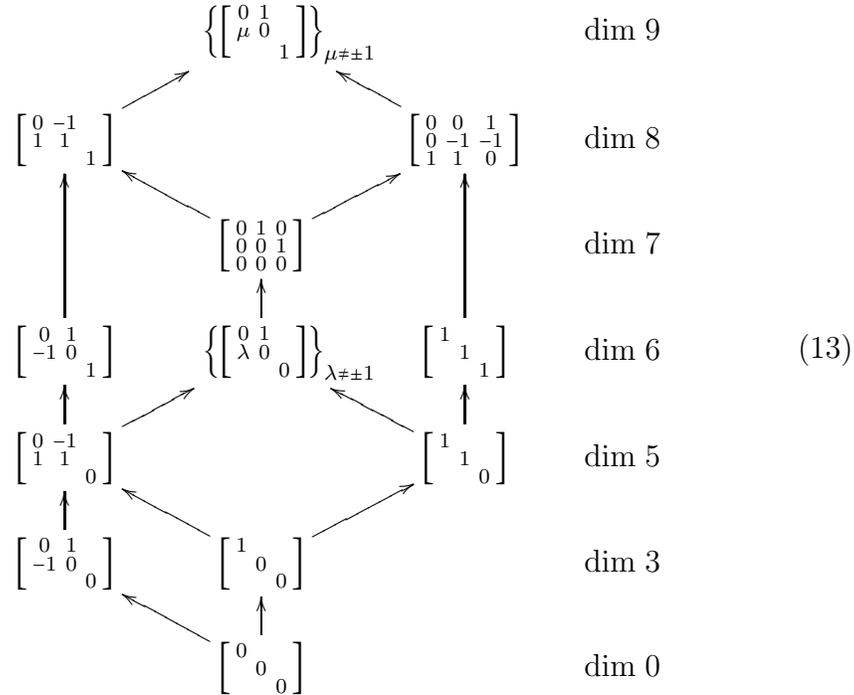
\begin{itemize}
  \item The bundle that
      corresponds to
      \[\left\{\left[\begin{smallmatrix}
      0&1&\\ {\lambda} &0&\\&&0
\end{smallmatrix}\right]
\right\}_{\lambda \ne\pm 1}\]
consists of all matrices whose
congruence canonical forms are
$\left[\begin{smallmatrix} 0&1&\\
{\lambda} &0&\\&&0
\end{smallmatrix}\right]$ with $\lambda \ne
\pm 1$.

  \item The bundle that
      corresponds to
      \[\left\{\left[\begin{smallmatrix}
      0&1&\\ {\mu} &0&\\&&1
\end{smallmatrix}\right]
\right\}_{\mu \ne\pm 1}\]
consists of all matrices whose
congruence canonical forms are
$\left[\begin{smallmatrix} 0&1&\\
{\mu} &0&\\&&1
\end{smallmatrix}\right]$
with $\mu \ne \pm 1$.

  \item The other bundles
      coincide with the
      congruence classes of the
      corresponding canonical
      matrices.
\end{itemize}
\end{itemize}
\end{theorem}

\begin{remark}\label{lke}
Since the number of congruence bundles
of matrices of fixed size $n\times n$
is finite, for each $A\in\mathbb
C^{n\times n}$ there exist congruence
bundles $\mathcal A_1,\dots,\mathcal
A_t\subset \mathbb C^{n\times n}$ such
that each sufficiently small
neighborhood of $A$ is contained in
${\mathcal A_1\cup\dots\cup\mathcal
A_t}$ and has a nonempty intersection
with each $\mathcal A_i$. If the
closure graph for congruence bundles of
$n\times n$ matrices is known, then the
set $\mathcal A_1,\dots,\mathcal A_t$
consists of all congruence bundles
$\mathcal A$ such that there exists a
directed path (including the ``lazy''
path of length $0$) from the vertex
that represents the bundle with $A$ to
the vertex that represents $\mathcal
A$. This information is important if
$A$ is known only approximately.

For example, the closure graph
$G_3^{\rm{B}}$ shows that each
sufficiently small neighborhood of
$I_3$ contains matrices whose
congruence canonical forms are
\[
\begin{bmatrix}
1&0&0\\ 0&1&0\\0&0&1
 \end{bmatrix},\quad
\begin{bmatrix}
0&0&1\\ 0&-1&-1\\1&1&0
 \end{bmatrix},\quad\text{and }
\begin{bmatrix}
0&1&0\\ \mu &0&0\\0&0&1
 \end{bmatrix}
\]
for some $\mu \ne\pm 1$, and does not
contain matrices with any other
congruence canonical form.

Arnold \cite[\S\,30E]{arn1} pointed to
possible applications of the partitions
into bundles: if in the study of a
phenomenon we obtain another partition,
then in the idealization of the
phenomenon something essential was
missed, or there were some special
reasons for an additional multiplicity
of parameters.

\end{remark}

\section{Supporting theory}
\label{ssk}

In this section we collect several
lemmas, which are used in the proofs of
the theorems.

\subsection{Roiter's assertion}

Instead of matrix pencils $A+\lambda
B$, it is convenient to consider matrix
pairs $(A,B)$ in which both matrices
have the same size. Matrix pairs
$(A,B)$ and $(A',B')$ are
\emph{equivalent} if there are
nonsingular $R$ and $S$ such that
$RAS=A'$ and $RBS=B'$ \cite{gan} (which
corresponds to the strict equivalence
of $A+\lambda B$ and $A'+\lambda B'$).
For each nonsingular matrix $A$,
$A^{-T}A:=(A^{-1})^TA$ is the
\emph{cosquare} of $A$.

The statement (a) in the following
lemma was formulated by Roiter
\cite{roi} and was extended to
arbitrary systems of linear mappings
and bilinear forms in
\cite{roi,ser_izv}. It reduces the
problem of classifying matrices up to
congruence to the problem of
classifying matrix pairs of the form
$(A,A^T)$ up to equivalence.

\begin{lemma}
\label{leq}
\begin{itemize}

  \item[{\rm(a)}] Two square
      matrices $A$ and $B$ are
      congruent if and only if the
      matrix pairs $(A,A^T)$ and
      $(B,B^T)$ are equivalent.

  \item[{\rm(b)}] Two nonsingular
      matrices $A$ and $B$ are
      congruent if and only if
      their cosquares $A^{-T}A$ and
      $B^{-T}B$ are similar.

\end{itemize}
\end{lemma}

\begin{proof} (a)
Let $A$ and $B$ be $n\times n$
matrices. Write $A=A_s+A_c$ and
$B=B_s+B_c$, where $A_s,B_s$ are
symmetric and $A_c,B_c$ are
skew-symmetric. Then there are three
equivalent conditions to state that $A$
and $B$ are congruent:
\begin{itemize}
  \item $(A_s,A_c)$ and $(B_s,B_c)$
      are congruent;
  \item $(A_s,A_c)$ and $(B_s,B_c)$
      are equivalent (see
      \cite[\S\,95, Theorem 3]{mal}
      or \cite[Lemma
      2.2]{f-h-s_trid});
  \item $(A, A^T)$ and $(B,B^T)$
      are equivalent.
\end{itemize}

(b) This statement follows from (a); it
is also proved in \cite[Lemma
2.1]{hor-ser_can} (or see \cite[Theorem
4.5.27]{hor-jo}).
\end{proof}

Note that Lemma \ref{leq}(a) is also
easily derived from the congruence
canonical form given in Proposition
\ref{le1}. Conversely, Proposition
\ref{le1} follows straightforwardly
from Lemma \ref{leq}(a) and the
Kronecker canonical form of matrix
pencils.

\subsection{Miniversal
deformations}\label{hyo}

We use a \emph{miniversal deformation}
of a square matrix $A$ under
congruence; that is, a normal form with
the minimal number of independent
parameters to which all matrices $A+X$
close to $A$ can be reduced by
transformations
\begin{equation*}\label{dtge}
A+X\mapsto {\cal
S}(X)^T (A+X) {\cal
S}(X),\qquad S(0)=I,
\end{equation*}
in which $S(X)$ depends holomorphically
on the entries of $X$. The notion of a
miniversal deformation was introduced
by Arnold \cite{arn}; he constructed
miniversal deformations of Jordan
canonical matrices.

A miniversal deformation of
$A\in{\mathbb C}^{n\times n}$ under
congruence was constructed in
\cite[Theorem 2.1]{f-ser} as follows.
Let $\cal D$ be any matrix with entries
$0$ and $*$ that satisfies the
condition
\begin{equation*}\label{jyr}
{\mathbb C}^{\,n\times
n}=T(A) \oplus {\cal
D}({\mathbb C}),
\end{equation*}
in which $T(A)$ is defined in
\eqref{msiy} and ${\cal D}({\mathbb
C})$ is the vector space of all
matrices obtained from $\cal D$ by
replacing its star entries by complex
numbers. Then all matrices in $A+{\cal
D}({\mathbb C})$ that are sufficiently
close to $A$ form a miniversal
deformation of $A$. The number
$n^2-\dim T(A)=\dim{\cal D}({\mathbb
C})$ is the {codimension of the
congruence class of $A$}; it is equal
to the number of stars in $\cal D$.

By \cite[Part III, Theorem 1.7]{tan},
the boundary of each congruence class
is a union of congruence classes of
strictly lower dimension, which ensures
the following lemma.

\begin{lemma}
\label{lej} If $v\to w$ is an arrow in
the closure graph $G_n$, then the
congruence class $C_v$ assigned to $v$
is contained in the closure of the
congruence class $C_w$, and so the
dimension of $C_v$ is lower than the
dimension of $C_w$.
\end{lemma}

The miniversal deformations of $2\times
2$ and $3\times 3$ matrices under
congruence are given in the following
lemma.

\begin{lemma}[{\cite[Example 2.1]{f-ser}}]
\label{le_def_a} Let $A$ be any
$2\times 2$ or $3\times 3$ matrix. Then
all matrices $A+X$ that are
sufficiently close to $A$ can be
simultaneously reduced by some
transformation
\begin{equation}\label{dtg}
{\cal
S}(X)^T (A+X) {\cal
S}(X),\quad\text{${\cal S}(X)$ is nonsingular and holomorphic
at $0$,}
\end{equation}
to one of the following forms, in which
$\lambda \in \mathbb
C\smallsetminus\{-1,1\}$ and each
nonzero $\lambda $ is determined up to
replacement by $\lambda^{-1} $:

$\rm{(i)}$ if $A$ is $2\times 2$,
\begin{align*}
&\begin{bmatrix} 0&\\
&0
 \end{bmatrix}+
\begin{bmatrix}
 *&*\\ *&*
 \end{bmatrix},
                            &&
\begin{bmatrix} 1&\\
&0
 \end{bmatrix}+
\begin{bmatrix}
 0&0\\ *&*
 \end{bmatrix},
                         &&
\begin{bmatrix} 1&\\
&1
 \end{bmatrix}+
\begin{bmatrix}
 0&0\\ *&0
 \end{bmatrix},
                          \\&
\begin{bmatrix} 0&1\\
-1&0
 \end{bmatrix}+
\begin{bmatrix}
 *&0\\ *&*
 \end{bmatrix},
                       &&
\begin{bmatrix} 0&-1\\
1&1
 \end{bmatrix}+
\begin{bmatrix}
 *&0\\ 0&0
 \end{bmatrix},
                         &&
\begin{bmatrix} 0&1\\
\lambda &0
 \end{bmatrix}+
\begin{bmatrix}
 0&0\\ *&0
 \end{bmatrix};
\end{align*}

$\rm{(ii)}$ if $A$ is $3\times 3$,
\begin{longtable}{ll}
$\begin{bmatrix}
0&&\\
&0&\\&&0
 \end{bmatrix}+
\begin{bmatrix}
 *&*&*\\  *&*&*\\ *&*&*
 \end{bmatrix},$
                             &
$\begin{bmatrix} 1&&\\
&0&\\&&0
 \end{bmatrix}+
\begin{bmatrix}
 0&0&0\\  *&*&*\\ *&*&*
 \end{bmatrix},$\vspace{3pt}
                           \\
$\begin{bmatrix} 1&&\\
&1&\\&&0
 \end{bmatrix}+
\begin{bmatrix}
 0&0&0\\  *&0&0\\ *&*&*
 \end{bmatrix},$
                            &
$\begin{bmatrix} 1&&\\ &1&\\&&1
 \end{bmatrix}+
\begin{bmatrix}
 0&0&0\\  *&0&0\\
 *&*&0
 \end{bmatrix},$\vspace{3pt}
                       \\
$\begin{bmatrix} 0&1&\\
-1&0&\\&&0
 \end{bmatrix}+
\begin{bmatrix}
 *&0&0\\  *&*&0\\ *&*&*
 \end{bmatrix},$
                      &
$\begin{bmatrix} 0&1&\\ \lambda
&0&\\&&0
 \end{bmatrix}+
\begin{bmatrix}
 0&0&0\\  *&0&0\\
 *&*&*
 \end{bmatrix}
\ (\lambda\ne 0),$\vspace{3pt}
                        \\
$\begin{bmatrix} 0&1&\\ 0 &0&\\&&0
 \end{bmatrix}+
\begin{bmatrix}
 0&0&0\\  *&0&*\\
 *&0&*
 \end{bmatrix},$
                          &
$\begin{bmatrix} 0&-1&\\ 1&1&\\&&0
 \end{bmatrix}+
\begin{bmatrix}
 *&0&0\\  0&0&0\\ *&*&*
 \end{bmatrix},$\vspace{3pt}
                           \\
$\begin{bmatrix} 0&1&\\
-1&0&\\&&1
 \end{bmatrix}+
\begin{bmatrix}
 *&0&0\\  *&*&0\\
 0&0&0
 \end{bmatrix},$
                         &
$\begin{bmatrix} 0&1&\\ \lambda
&0&\\&&1
 \end{bmatrix}+
\begin{bmatrix}
 0&0&0\\  *&0&0\\
 0&0&0
 \end{bmatrix},$\vspace{3pt}
                   \\
$ \begin{bmatrix} 0&-1&\\ 1&1&\\&&1
 \end{bmatrix}+
\begin{bmatrix}
 *&0&0\\  0&0&0\\
 0&0&0
 \end{bmatrix},$
                      &
$\begin{bmatrix} 0&1&0\\ 0&0&1\\0&0&0
 \end{bmatrix}+
\begin{bmatrix}
 0&0&0\\  0&0&0\\
 *&0&*
 \end{bmatrix},$\vspace{3pt}
                     \\
$ \begin{bmatrix} 0&0&1\\
0&-1&-1\\1&1&0
 \end{bmatrix}+
\begin{bmatrix}
 0&0&0\\  *&0&0\\
 0&0&0
 \end{bmatrix}.$
\end{longtable}

Each of these matrices has the form
$A_{\rm can}+{\cal D}$ in which $A_{\rm
can}$ is a canonical matrix for
congruence and the stars in ${\cal D}$
are complex numbers that tend to zero
as $X$ tends to zero. The number of
stars is the smallest that can be
attained by using transformations
\eqref{dtg}; it is equal to the
codimension of the congruence class of
$A$.
\end{lemma}

\subsection{Characteristic polynomial}
\label{s_pr0}

\begin{lemma}\label{lemz}
Let $A$ and $B$ be nonsingular
matrices. If there is an arbitrarily
small perturbation of $A$ that is
congruent to $B$ then the cosquares
$A^{-T}A$ and $B^{-T}B$ have the same
characteristic polynomial and so
$\det(A^Tx-A)=c\cdot\det(B^Tx-B)$ for a
nonzero $c\in\mathbb C$.
\end{lemma}

\begin{proof}
There exists a sequence of nonsingular
matrices $A_1,\,A_2,\,\dots$ converging
to $A$, in which all $A_i$ are
congruent to $B$. By Proposition
\ref{leq}(b), all $A_i^{-T}A_i$ are
similar to $B^{-T}B$ and so they have
the same characteristic polynomial
$\chi (x)$. Since
$A_1^{-T}A_1,\,A_2^{-T}A_2,\,\dots$
converges to $A^{-T}A$, $\chi (x)$ is
the characteristic polynomial of
$A^{-T}A$.
\end{proof}

\section{Proof
 of the theorems}
\label{sss}

\subsection{Proof of Theorem \ref{the1}}
\label{s_pr1}

For each $2\times 2$ congruence
canonical matrix $M,$ the dimension
$d_M$ of its congruence class is
indicated in \eqref{g1}. It was
calculated as follows: \eqref{kow}
ensures that $d_M=4-c_M$ in which $c_M$
is the codimension of the congruence
class of $M$; $c_M$ is equal to the
number of stars in $M+{\cal D}$ from
Lemma \ref{le_def_a}(i).

The proof of Theorem \ref{the1} is
divided into two steps.

\emph{Step 1:} Let us prove that each
arrow $M\to N$ in \eqref{g1} is
correct; that is, there is an
arbitrarily small perturbation of the
matrix $M$ that is congruent to $N$.

We write $A\cong  B$ if $A$ is
congruent to $B$; that is, $S^TAS=B$
for some nonsingular $S$. Expressing
$S$ as a product of elementary
matrices, we obtain that
\begin{equation}\label{jyr1}
\parbox[c]{0.7\textwidth}{$A\cong  B$ if and
only if $A$ can be
reduced to $B$ by a sequence of
elementary
transformations of
rows and the same
transformations of columns.
}\end{equation}

We denote by $\varepsilon,\ \delta$,
and $\zeta$ arbitrarily small complex
numbers.

$\bullet$ ${\left[\begin{smallmatrix}
0&0\\
0&0
 \end{smallmatrix}\right]}\to{\left[\begin{smallmatrix}
0&1\\
-1&0
 \end{smallmatrix}\right]} $ and
${\left[\begin{smallmatrix}
0&0\\
0&0
 \end{smallmatrix}\right]}\to{\left[\begin{smallmatrix} 1&0\\
0&0
 \end{smallmatrix}\right]}$ since
\[
\begin{bmatrix} 0&\varepsilon
\\
-\varepsilon &0
 \end{bmatrix}\cong
 \begin{bmatrix} 0&1\\
-1 &0
 \end{bmatrix}
 \quad\text{and}\quad
 \begin{bmatrix}
\varepsilon &0\\0&0
 \end{bmatrix}\cong
 \begin{bmatrix} 1&0\\
0&0
 \end{bmatrix}.
\]

$\bullet$  ${\left[\begin{smallmatrix}
0&1\\
-1&0
 \end{smallmatrix}\right]}
\to{\left[\begin{smallmatrix} 0&-1\\
1&1
 \end{smallmatrix}\right]}
$ since
\[
\begin{bmatrix}
0&1
\\
-1 &\varepsilon
 \end{bmatrix}\cong
 \begin{bmatrix} 0&-1\\
1 &1
 \end{bmatrix}.
\]

$\bullet$ ${\left[\begin{smallmatrix} 1&0\\
0&0
 \end{smallmatrix}\right]}\to {\left[\begin{smallmatrix} 0&-1\\ 1&1
 \end{smallmatrix}\right]}
$ since
\[
 \begin{bmatrix}
1 &\varepsilon\\
-\varepsilon&0
 \end{bmatrix}\cong
 \begin{bmatrix} 1&1\\
-1&0
 \end{bmatrix} \quad\text{and}\quad
\begin{bmatrix} 0&1\\
1&0
 \end{bmatrix}\begin{bmatrix} 1&1\\
-1&0
 \end{bmatrix}\begin{bmatrix} 0&1\\
1&0
 \end{bmatrix}=\begin{bmatrix} 0&-1\\
1&1
 \end{bmatrix}.
\]

$\bullet$  ${\left[\begin{smallmatrix} 1&0\\
0&0
 \end{smallmatrix}\right]}\to
\left[\begin{smallmatrix}
0&1\\
0&0
 \end{smallmatrix}\right]
$
and ${\left[\begin{smallmatrix} 1&0\\
0&0
 \end{smallmatrix}\right]}\to {I_2}$ since
\[
 \begin{bmatrix}
1&\varepsilon
\\
0&0
 \end{bmatrix}\cong
 \begin{bmatrix} 0&1\\
0&0
 \end{bmatrix}
 \quad\text{and}\quad
\begin{bmatrix}
1 &0\\
0&\varepsilon
 \end{bmatrix}\cong
 \begin{bmatrix} 1&0\\
0&1
 \end{bmatrix}.
\]

$\bullet$ ${\left[\begin{smallmatrix} 1&0\\
0&0
 \end{smallmatrix}\right]}\to
\left[\begin{smallmatrix}
0&1\\
\lambda &0
 \end{smallmatrix}\right]
$ for each $\lambda\notin \{-1,0, 1\}$;
it suffices to find arbitrarily small
$\delta\ne 0 $ and $\varepsilon$ for
which
\begin{equation}\label{htx}
A :=
 \begin{bmatrix}
1 &0\\
\varepsilon&\delta
 \end{bmatrix}\quad \text{is congruent to}\quad
H_2(\lambda )=
 \begin{bmatrix}
0 &1\\
\lambda &0
 \end{bmatrix},
\end{equation}
i.e. $A^{-T}A$ is similar to
$H_2(\lambda )^{-T}H_2(\lambda
)=\diag(\lambda, \lambda^{-1})$ (see
Proposition \ref{leq}(b)), i.e. the
eigenvalues of $A^{-T}A$ are $\lambda$
and $\lambda ^{-1}$, i.e. the roots of
the polynomial
\[
\det(A^Tx-A)=
\begin{vmatrix}
x-1 &\varepsilon x\\
-\varepsilon&\delta x-
\delta
 \end{vmatrix}=\delta
 [x^2+(\varepsilon
 ^2/\delta-2
 )x+1]
\]
are $\lambda$ and $\lambda ^{-1}$. This
polynomial is equal to $\delta
(x-\lambda)( x-\lambda ^{-1})$ if
\begin{equation}\label{ldt}
\varepsilon^2/\delta-2=
-\lambda-\lambda ^{-1},
\end{equation}
which can be satisfied by some
arbitrarily small $\varepsilon $ and
$\delta $.

\emph{Step 2:} Let us prove that we
have not missed any arrows in
\eqref{g1}.

We write $M\nrightarrow N$ if the
closure graph does not have the arrow
$M\to N$; i.e., if each matrix obtained
from $M$ by an arbitrarily small
perturbation is not congruent to $N$.
The evident statement ``if $M
\rightarrow N$ and $M\nrightarrow P$,
then $N\nrightarrow P$'' and Lemma
\ref{lej} ensure that we need to prove
only the absence of the arrows $H_2(-1)
\to H_2(\lambda)$ and $H_2(-1) \to
I_2$.

$\bullet$ ${H_2(-1)} \nrightarrow
H_2(0)$ since each arbitrarily small
perturbation of $ H_2(-1) $ has rank
$2$ and $ H_2(0) $ has rank $1$.

$\bullet$ $H_2(-1) \nrightarrow
{H_2}(\lambda)$ ($\lambda \ne 0,1,-1$)
by Lemma \ref{lemz} since
$\det(H_2(-1)^Tx-H_2(-1))=(x+1)^2$ and
$\det(H_2(\lambda)^{T}x-H_2(\lambda))=
-(x-\lambda )(\lambda x-1)$.

$\bullet$ $H_2(-1)\nrightarrow I_2$.
Each arbitrarily small perturbation of
$H_2(-1) $ is not congruent to $I_2$
since their skew-symmetric parts are
not congruent.\hfill\hbox{\qedsymbol}

\subsection{Proof of Theorem \ref{the2}}
\label{s_pr2}

Each canonical matrix for congruence is
a direct sum of blocks of the forms
$H_{2m}(\lambda)$, $\Gamma_n$, and
$J_k(0)$ (see \eqref{can}). Replacing
them by $[\lambda]^{2m}$, $1^n$, and
$0^k$, respectively, and deleting the
symbols $\oplus$, we get a compact
notation of the canonical matrix. For
example,
\[
[1]^4[3]^21^61^610^20\quad\text{is}\quad
H_4(1)\oplus H_2(3)\oplus
\Gamma_6\oplus \Gamma_6\oplus \Gamma_1
\oplus J_2(0)\oplus J_1(0)
\]
(we write $1$ and $0$ instead of $1^1$
and $0^1$).

All $3\times 3$ canonical matrices for
congruence are written in this notation
as follows:
\begin{equation}\label{g2v}
\begin{split}
\begin{aligned}
000=&\ \begin{bmatrix}
0&&\\ &0&\\&&0
 \end{bmatrix},
  & [-1]^20=&\
\begin{bmatrix}
0&1&\\ -1&0&\\&&0
 \end{bmatrix},
&100=&\ \begin{bmatrix}
1&&\\ &0&\\&&0
 \end{bmatrix},\\
  1^20=&\
\begin{bmatrix}
0&-1&\\ 1&1&\\&&0
 \end{bmatrix},
           &
[\lambda]^20=&\
\begin{bmatrix}
0&1&\\ \lambda
&0&\\&&0
 \end{bmatrix},
  & 110=&\
\begin{bmatrix}
1&&\\ &1&\\&&0
 \end{bmatrix},\\
[-1]^21=&\ \begin{bmatrix}
0&1&\\ -1&0&\\&&1
 \end{bmatrix},
  & 111=&\
\begin{bmatrix}
1&&\\ &1&\\&&1
 \end{bmatrix},
 &0^3=&\ \begin{bmatrix}
0&1&0\\ 0&0&1\\0&0&0
 \end{bmatrix},\\
  1^21=&\
\begin{bmatrix}
0&-1&\\ 1&1&\\&&1
 \end{bmatrix},
 &
[\mu]^21=&\
\begin{bmatrix}
0&1&\\ \mu &0&\\&&1
 \end{bmatrix},
\ & 1^3=&\
\begin{bmatrix}
0&0&1\\ 0&-1&-1\\1&1&0
 \end{bmatrix}.
 \end{aligned}
\end{split}
\end{equation}
They are direct sums of blocks of the
form \eqref{can};  the zero entries
outside of these blocks are not shown.

For each $3\times 3$ canonical matrix
$M$, the dimension $d_M$ of its
congruence class is indicated in
\eqref{g3}. It was calculated as
follows: \eqref{kow} ensures that
$d_M=9-c_M$ in which $c_M$ is the
codimension of the congruence class of
$M$; $c_M$ is equal to the number of
stars in $M+{\cal D}$ from Lemma
\ref{le_def_a}(ii).

\emph{Step 1:} Let us prove that each
arrow in \eqref{g3} is correct.

For any ${2\times 2}$ matrices $A$ and
$B$ and any $c\in\mathbb C$, if each
neighborhood of $A$ contains a matrix
that is congruent to $B$, then each
neighborhood of $A\oplus [c]$ contains
a matrix that is congruent to $B\oplus
[c]$. Thus, the arrows of \eqref{g1}
ensure the correctness of the following
arrows in \eqref{g3}:
\begin{align*}
 000\to [-1]^20&\text{ since }
{\left[\begin{smallmatrix}
0&0\\
0&0
 \end{smallmatrix}\right]}\to{\left[\begin{smallmatrix} 0&1\\
-1&0
 \end{smallmatrix}\right]}
& 100\to 110&\text{ since
}{\left[\begin{smallmatrix} 1&0\\
0&0
 \end{smallmatrix}\right]}\to{I_2}
       \\
 000\to 100&\text{ since }
{\left[\begin{smallmatrix}
0&0\\
0&0
 \end{smallmatrix}\right]}\to{\left[\begin{smallmatrix} 1&0\\
0&0
 \end{smallmatrix}\right]}
&
 [-1]^20\to 1^20&\text{ since }{
  \left[\begin{smallmatrix} 0&1\\
-1&0
 \end{smallmatrix}\right]}
\to{\left[\begin{smallmatrix} 0&-1\\ 1&1
 \end{smallmatrix}\right]}
      \\
110\to 111&\text{ since }{
\left[\begin{smallmatrix} 1&0\\
0&0
 \end{smallmatrix}\right]}\to{I_2}
        &
 100\to 1^20&\text{ since }{\left[\begin{smallmatrix} 1&0\\
0&0
 \end{smallmatrix}\right]}\to{\left[\begin{smallmatrix} 0&-1\\ 1&1
 \end{smallmatrix}\right]}
         \\
[-1]^21\to
1^21&\text{ since }{\left[\begin{smallmatrix} 0&1\\
-1&0
 \end{smallmatrix}\right]}
\to{\left[\begin{smallmatrix} 0&-1\\ 1&1
 \end{smallmatrix}\right]}
         &
100\to [\lambda]^20
&\text{ since }{\left[\begin{smallmatrix} 1&0\\
0&0
 \end{smallmatrix}\right]}\to \left[\begin{smallmatrix} 0&1\\
\lambda &0
 \end{smallmatrix}\right].
\end{align*}
The correctness of the remaining arrows
in \eqref{g3} is proved as follows.

$\bullet$ $1^20\to [-1]^21$ since
\[
\begin{bmatrix}
-1&0&0
\\
 0&1&-1/\varepsilon
\\ 1&0&1/\varepsilon
 \end{bmatrix}
\begin{bmatrix}
0&-1&0\\1&1&0\\0&2\varepsilon
&\varepsilon^2 \end{bmatrix}
\begin{bmatrix}
-1&0&1\\0&1&0
\\ 0&-1/\varepsilon&1/\varepsilon
 \end{bmatrix}
=
\begin{bmatrix}
0&1&0\\-1&0&0\\0&0&1
 \end{bmatrix}.
\]

$\bullet$ $1^20\to 0^3$ since we can
reduce
\begin{equation}\label{kic}
\begin{bmatrix}
0&-1&0
\\
 1&1 &\varepsilon
\\ 0&0&0
 \end{bmatrix}
\quad\text{to}\quad
\begin{bmatrix}
0&-1&0
\\
0 &0& \varepsilon  \\
0&0&0
 \end{bmatrix}
\end{equation}
by transformations \eqref{jyr1}: we add
the last column multiplied by
$-\varepsilon ^{-1}$ to the first and
second columns; the same
transformations of rows do not change
the matrix.

$\bullet$ $\{[\lambda]^20,110\}\to 0^3$
(i.e., $[\lambda]^20\to 0^3$ and
$110\to 0^3$) since
\[
\begin{bmatrix}
0&1&0
\\
\lambda &0 &\varepsilon \\ 0&0&0
 \end{bmatrix}\cong
\begin{bmatrix}
0&1&0
\\
0&0& \varepsilon  \\
0&0&0
 \end{bmatrix}\]
 and
\[\begin{bmatrix}
1&i&0
\\
0 &1&0 \\
0&0&1
 \end{bmatrix}
\begin{bmatrix}
1&0&-i\varepsilon
\\
0&1&\varepsilon  \\
0&0&0
 \end{bmatrix}
 \begin{bmatrix}
1&0&0
\\
i &1&0 \\
0&0&1
 \end{bmatrix}
=\begin{bmatrix} 0&i&0
\\
i &1&\varepsilon  \\
0&0&0
 \end{bmatrix}
 \cong
\begin{bmatrix}
0&i&0
\\
0&0& \varepsilon  \\
0&0&0
 \end{bmatrix},
\]
as in \eqref{kic}.

$\bullet$ $111\to 1^3$ since
\[
\begin{bmatrix}
1 &0&0
\\
0&1&0  \\
0&0&1
 \end{bmatrix}\cong
\begin{bmatrix}
0&0&1
\\
0&-1& 0 \\
1&0&0
 \end{bmatrix}\]
 by Proposition
\ref{leq}(b),  and
\[
\begin{bmatrix}
\varepsilon &0&0
\\
0&1&0  \\
0&0&\varepsilon^{-1}
 \end{bmatrix}
\begin{bmatrix}
0&0&1
\\
0&-1&-\varepsilon  \\
1&\varepsilon&0
 \end{bmatrix}
 \begin{bmatrix}
\varepsilon &0&0
\\
0&1&0  \\
0&0&\varepsilon^{-1}
\end{bmatrix}
=\begin{bmatrix} 0&0&1
\\
0&-1 &-1  \\
1&1&0
 \end{bmatrix}.
\]
$\bullet$ $0^3\to
\{1^21,[\mu]^21,1^3\}$. For the
perturbation
\[
A:=\begin{bmatrix}
0&1&0
\\
0 &0&1  \\
\varepsilon &0&\delta
 \end{bmatrix}\] of $0^3$, we have
\begin{equation} \label{gpe}
\begin{split}
\det(A^Tx-A)&=\begin{vmatrix}
0&-1&\varepsilon x
\\
x &0&-1  \\
-\varepsilon &x&\delta
x-\delta
 \end{vmatrix}
 =\varepsilon x^3-\varepsilon
+x(\delta x-\delta)\\
&=\varepsilon
(x^3+\delta
\varepsilon^{-1}x^2-
\delta
\varepsilon^{-1}x-1)\\
&=\varepsilon(x-1)(x^2+(\delta
\varepsilon^{-1}+1)x+1).
\end{split}
\end{equation}

To verify $0^3\to 1^21$, we set $\delta
=\varepsilon $ and obtain $
\det(A^Tx-A)=\varepsilon (x-1)(x+1)^2.
$ Thus, the eigenvalues of $A^{-T}A$
are $1,-1,-1$. Only the cosquares of
$[-1]^21$ and $1^21$ among the
cosquares of the nonsingular matrices
in \eqref{g2v} have these eigenvalues;
in particular, because the cosquare of
$1^3$ is
\begin{equation}\label{kef}
\Gamma _3^{-T}\Gamma _3=\begin{bmatrix}
                    1&-1 & 1 \\
                    1 & -1&0 \\1&0&0
                  \end{bmatrix}
\begin{bmatrix}
                    0&0 & 1 \\
                    0& -1&-1 \\1&1&0
                  \end{bmatrix}
=\begin{bmatrix}
                    1&2 & 2 \\
                    0&1&2 \\0&0&1
                  \end{bmatrix}.
\end{equation}
Proposition \ref{leq}(b) ensures that
$A$ is congruent to one of these
matrices; that is, $0^3\to [-1]^21$ or
$0^3\to 1^21$. But $0^3\nrightarrow
[-1]^21$ by Lemma \ref{lej}.

Let us verify $0^3\to [\mu]^21$. Each
complex number is represented in the
form $\delta \varepsilon^{-1}$ with
arbitrarily small $\delta$ and
$\varepsilon$; hence for every nonzero
$\mu $ there exist arbitrarily small
$\delta$ and $\varepsilon$ such that $
x^2+(\delta \varepsilon^{-1}+1)x+1=
(x-\mu )(x-\mu^{-1}). $ By \eqref{gpe},
the eigenvalues of $A^{-T}A$ are $1,\mu
,\mu ^{-1}$. In the list \eqref{g2v},
only the cosquare of $[\mu]^21$ has
these eigenvalues. Proposition
\ref{leq}(b) ensures that $0^3\to
[\mu]^21$ if $\mu \ne 0$. The arrow
$0^3\to [0]^21$ exists because
\[
\begin{bmatrix}
0&1&0
\\
0 &0&1  \\
0&0&\delta
 \end{bmatrix}
\cong
\begin{bmatrix}
0&1&0
\\
0 &0&0 \\
0 &-1&\delta
 \end{bmatrix}
\cong
\begin{bmatrix}
0&1&0
\\
0 &0&0 \\
0 &0&\delta
 \end{bmatrix}
 \cong
\begin{bmatrix}
0&1&0
\\
0 &0&0 \\
0 &0&1
 \end{bmatrix}.
\]

To verify $0^3\to 1^3$, we set $\delta
=-3\varepsilon $. Then
$\det(A^Tx-A)=\varepsilon (x-1)^3$. By
\eqref{kef} and Proposition
\ref{leq}(b), $A$ is congruent to $111$
or $1^3$ from \eqref{g2v}. But
$0^3\nrightarrow 111$ by Lemma
\ref{lej}.

\emph{Step 2:} Let us prove that we
have not missed any arrows in
\eqref{g3}.

Due to Lemma \ref{lej}, it suffices to
verify that
\begin{align*}
[-1]^20&\nrightarrow
\{[\lambda]^20,111\}&
\{[\lambda]^20,110\} &\nrightarrow
[-1]^21\\
[\lambda]^20
&\nrightarrow 111 &[-1]^21&\nrightarrow
\{[\mu]^21,1^3\}\\
111&\nrightarrow\{1^21,[\mu]^21\}
\end{align*}
for all $\lambda \ne \pm 1$ and $\mu
\ne \pm 1$.

$\bullet$ $[-1]^20\nrightarrow
[\lambda]^20$ $(\lambda \ne \pm 1)$. To
the contrary, suppose that there is an
arbitrarily small perturbation $A$ of
$[-1]^20$ that is congruent to
$[\lambda]^20$ for some $\lambda \ne
\pm 1$. We can suppose that this
perturbation has the form given in
Lemma \ref{le_def_a}(ii); that is,
\begin{equation}\label{ltd}
A:=\begin{bmatrix}
\alpha &1&0\\
-1+\beta  &\gamma &0\\
\varepsilon &\zeta&
\eta
 \end{bmatrix}
\cong  \begin{bmatrix}
0 &1&0\\
\lambda   &0 &0\\
0 &0 &0
 \end{bmatrix}
\end{equation}
for some arbitrarily small
$\alpha,\beta,\gamma,\varepsilon,
\zeta, \eta$. Then $\rank A=2$, and so
$\eta=0$.

If $\varepsilon = \zeta=0$, then the
congruence \eqref{ltd} with $\eta=0$
implies
\[
\begin{bmatrix}
\alpha &1\\
-1+\beta  &\gamma
 \end{bmatrix}
\cong  \begin{bmatrix}
0 &1\\
\lambda   &0
 \end{bmatrix}
\]
because of the uniqueness in
Proposition \ref{le1}, which is
impossible since
$\left[\begin{smallmatrix} 0&1\\
-1&0
 \end{smallmatrix}\right] \nrightarrow
\left[\begin{smallmatrix} 0&1\\
\lambda &0
 \end{smallmatrix}\right]$
in \eqref{g1}.

Let $\zeta\ne 0$ and $|\varepsilon|\le
|\zeta|$.  By adding the second column
of $A$ to the first, we make
$\varepsilon=0$; then repeat this
transformation with rows. The entries
$\alpha $ and $\beta $ are changed to
$\alpha'$ and $\beta'$, which remain
arbitrarily small. By adding the last
row, we make the entries above $\zeta$
equal to zero. Thus,
\begin{equation*}\label{dttj}
A\cong  \begin{bmatrix}
\alpha' &0&0\\
-1+\beta'  &0 &0\\
0 &\zeta& 0
 \end{bmatrix}
\cong  \begin{bmatrix}
\alpha' &0&0\\
1&0 &0\\
0 &1& 0
 \end{bmatrix}.
\end{equation*}
If $\alpha'=0$, then the last matrix is
congruent to $0^3$, which is not
congruent to $[\lambda]^20$. If
$\alpha'\ne 0$, then
\begin{equation*}\label{dtj}
A\cong  \begin{bmatrix}
\alpha' &*&0\\
0  &0 &0\\
0 &1& 0
 \end{bmatrix}
\cong  \begin{bmatrix}
\alpha' &0&0\\
0&0 &0\\
0 &1& 0
 \end{bmatrix}\cong  \begin{bmatrix}
1&0&0\\
0&0 &0\\
0 &1& 0
 \end{bmatrix},
\end{equation*}
which is congruent to $[0]^21$ and so
it is not congruent to $[\lambda]^20$.

Let $|\varepsilon|> |\zeta|$. We
interchange the first two rows and the
first two columns in $A$ and reason as
in the previous case.

$\bullet$ $[-1]^20\nrightarrow 111$
since if a perturbation of $[-1]^20$ is
congruent to $111$ then their
skew-symmetric parts must be congruent
too.

$\bullet$ $\{[\lambda]^20,110\}
\nrightarrow [-1]^21$; that is,
\[
\begin{bmatrix}
0&1&0\\ \lambda
&0&0\\0&0&0
 \end{bmatrix}+E
\not\cong
\begin{bmatrix}
0&1&0\\ -1&0&0\\0&0&1
 \end{bmatrix}
 \quad\text{and}\quad
\begin{bmatrix}
1&0&0\\ 0&1&0\\0&0&0
 \end{bmatrix}+E
\not\cong
\begin{bmatrix}
0&1&0\\ -1&0&0\\0&0&1
 \end{bmatrix}
\]
for $\lambda \ne \pm 1$ and each
sufficiently small matrix $E$. This is
true since the symmetric parts of the
corresponding matrices have unequal
ranks.

$\bullet$ $[\lambda]^20 \nrightarrow
111$; that is,
\[
\begin{bmatrix}
0&1&0\\ \lambda
&0&0\\0&0&0
\end{bmatrix}+E \not\cong
\begin{bmatrix}
1&0&0\\ 0&1&0\\0&0&1
 \end{bmatrix}
\]
for $\lambda \ne\pm 1$ and each
sufficiently small matrix $E$. This is
true since the skew-symmetric parts of
the corresponding matrices have unequal
ranks.

$\bullet$ $[-1]^21 \nrightarrow
\{[\mu]^21,1^3\}$; that is,
\[
\begin{bmatrix}
0&1&0\\ -1&0&0\\0&0&1
 \end{bmatrix}
+E
\not\cong
\begin{bmatrix}
0&1&0\\ \mu &0&0\\0&0&1
 \end{bmatrix}
       \quad\text{and}\quad
\begin{bmatrix}
0&1&0\\ -1&0&0\\0&0&1
 \end{bmatrix}
+E
 \not\cong
\begin{bmatrix}
0&0&1\\ 0&-1&-1\\1&1&0
 \end{bmatrix}
\]
for $\mu \ne\pm 1$ and each
sufficiently small matrix $E$. This is
true by Lemma \ref{lemz} since the
characteristic polynomials of the
cosquares of $[-1]^21$, $[\mu]^21$, and
$1^3$ are $(x+1)^2(x-1)$,
$(x-\mu)(x-\mu^{-1})(x-1)$, and
$(x-1)^3$ (see \eqref{kef}),
respectively.

$\bullet$
$111\nrightarrow\{1^21,[\mu]^21\}$;
that is,
\[
\begin{bmatrix}
1&0&0\\ 0&1&0\\0&0&1
 \end{bmatrix}
+E
\not\cong
\begin{bmatrix}
0&-1&0\\ 1 &1&0\\0&0&1
 \end{bmatrix}
       \quad\text{and}\quad
\begin{bmatrix}
1&0&0\\ 0&1&0\\0&0&1
 \end{bmatrix}
+E
 \not\cong
\begin{bmatrix}
0&1&0\\ \mu &0&0\\0&0&1
 \end{bmatrix}
 \]
for $\mu \ne\pm 1$ and each
sufficiently small matrix $E$. This is
true by Lemma \ref{lemz} since the
characteristic polynomials of the
cosquares of $111$, $1^21$, and
$[\mu]^21$ are $(x-1)^3$,
$(x+1)^2(x-1)$, and
$(x-\mu)(x-\mu^{-1})(x-1)$,
respectively.\hfill\hbox{\qedsymbol}

\subsection{Proof of Theorem
\ref{theorm}} \label{ghy} The Kronecker
canonical form for matrix pencils
ensures that each matrix pair is
equivalent to a direct sum, determined
uniquely up to permutation of summands,
of pairs of the form
\begin{equation*}\label{ksca}
(L_r,R_r),\quad
(L_r^T,R_r^T),\quad {P}_r(\lambda ):=
\begin{cases}
(I_r,J_r(\lambda)),&\text{if $\lambda \in
\mathbb C$},\\
(J_r(0),I_r),&\text{if $\lambda =
\infty,$}
\end{cases}
\end{equation*}
in which $L_r$ and $R_r$ are defined in
\eqref{1.4}, $r=1,2,\dots$, and
$\lambda \in\mathbb C\cup\infty$. This
sum is called the \emph{Kronecker
canonical form} of $(A,B)$ and is
denoted by $(A,B)\ca$.

The following definition in terms of
pencils was given by Edelman, Elmroth,
and K\r{a}gstr\"{o}m \cite[Section
3.1]{kag2}.

\begin{definition}\label{moy}
Let $(A,B)$ be any pair of matrices of
the same size. Write its Kronecker
canonical form as follows:
\begin{equation}\label{lui}
(A,B)\ca=
{\cal K}\oplus {\cal P}_1(\lambda_1)
\oplus\dots\oplus {\cal P}_t(\lambda_t),
\quad \lambda_i\ne \lambda _j\text{ if }
i\ne j,\quad \lambda_1,\dots,\lambda_t
\in\mathbb C\cup\infty,
\end{equation}
in which ${\cal K}$ is a direct sum of
pairs of the form $(L_k,R_k)$ and
$(L_k^T,R_k^T)$, and each ${\cal
P}_i(\lambda_i)$ is a direct sum of
pairs of the form ${P}_k(\lambda_i)$.
Then the \emph{equivalence bundle} of
$(A,B)$ consists of all matrix pairs
whose Kronecker canonical form is
\begin{equation}\label{mjx}
{\cal K}\oplus {\cal P}_1(\mu_1)
\oplus\dots\oplus {\cal P}_t(\mu_t),
\quad \mu_i\ne \mu_j\text{ if }
i\ne j,\quad \mu_1,\dots,\mu_t
\in\mathbb C\cup\infty
\end{equation}
with the same ${\cal K},{\cal P}_1,
\dots, {\cal P}_t.$
\end{definition}

Definition \ref{kuq1} of congruence
bundles can be reformulated as follows.

\begin{definition}\label{mos1}
The \emph{congruence bundle} of a
square matrix $A$ consists of all
matrices $B$ such that the pairs
$(A,A^T)$ and $(B,B^T)$ belong to the
same equivalence bundle.
\end{definition}

\begin{proof}[Proof of Theorem
\ref{theorm}] Let $B=B_1\oplus
\dots\oplus B_s$ be a direct sum of
matrices of the form $H_{2m}(\lambda),$
$\Gamma _n,$ and $J_{2r-1}(0)$; see
\eqref{cjy}. Then $(B,B^T)\ca$ can be
obtained by replacing in $B_1\oplus
\dots\oplus B_s$ each summand $B_i$ by
$(B_i,B_i^T)\ca$, which is calculated
as follows:
\begin{itemize}
  \item
      $(H_{2m}(\lambda),H_{2m}(\lambda)^T)
      \ca=P_m(\lambda)\oplus
      P_m(\lambda^{-1})$ (we set
      $0^{-1}:=\infty$);

  \item $(\Gamma _n,\Gamma
      _n^T)\ca=P_n((-1)^{n+1})$
      since $(\Gamma _n,\Gamma
      _n^T)$ is equivalent to
      $(\Gamma _n^{-T}\Gamma
      _n,I_n)$, which is equivalent
      to $(J_n((-1)^{n+1}),I_n)$
      because
\begin{equation*}\label{1nh}
\Gamma_n^{-1}=(-1)^{n+1}
 \begin{bmatrix}
\vdots&\vdots&\vdots&\vdots&\udots
\\
-1&-1&-1&-1&\\ 1&1&1&&\\ -1&-1&&&\\
1&&&& 0
\end{bmatrix},\quad
\Gamma_n^{-T}\Gamma_n=
(-1)^{n+1}
\begin{bmatrix} 1&2&&\text{\raisebox{-6pt}
{\large\rm *}}
\\&1&\ddots&\\
&&\ddots&2\\
0 &&&1
\end{bmatrix};
\end{equation*}

  \item $ (J_k(0),J_k(0)^T)\ca=
  \begin{cases}
(L_r,R_r)\oplus(L_r^T,R_r^T)
 & \hbox{if }k =2r-1 \\
 P_r(0)\oplus
P_r(\infty) & \hbox{if }k=2r
  \end{cases}
$\\ since $J_{k}(0)$ is
permutationally congruent to
\eqref{kuf}.
\end{itemize}

Let $B'$ be another direct sum of
matrices of the form \eqref{cjy}. Then
there are three equivalent conditions
to state that the matrices $B$ and $B'$
are in the same congruence bundle:
\begin{itemize}
  \item $(B,B^T)$ and $(B',{B'}^T)$
      are in the same equivalence
      bundle;
  \item $(B',{B'}^T)\ca$ can be
      obtained from $(B,{B}^T)\ca$
      as \eqref{mjx} from
      \eqref{lui};
  \item $B'$ can be obtained from
      $B$ as \eqref{kdr} from
      \eqref{lmio}.{\qedhere}
\end{itemize}
\end{proof}

\subsection{Proof of Theorem
\ref{the1_b}} \label{gep}

(a)  The vertices of graph \eqref{g2}
represent the congruence bundles of
$2\times 2$ matrices, which are
described in Theorem \ref{theorm}.
Deleting the
vertices $\left[\begin{smallmatrix} 0&1\\
\lambda &0
 \end{smallmatrix}\right]$
 ($\lambda \ne\pm 1$)
and the arrows
$\left[\begin{smallmatrix} 1&0\\
0 &0
 \end{smallmatrix}\right]
\to \left[\begin{smallmatrix} 0&1\\
\lambda &0
 \end{smallmatrix}\right]$ in the closure graph $G_2$
constructed in \eqref{g1}, we obtain a
subgraph of  \eqref{g2}. This proves
that all the arrows in \eqref{g2} are
correct with the possible exception of
\[{ \left[\begin{smallmatrix} 0&-1\\
1&1
 \end{smallmatrix}\right]}
\to\{\left[\begin{smallmatrix} 0&1\\
\lambda &0
 \end{smallmatrix}\right]\}_{\lambda
 \ne\pm 1}\quad \text{and}\quad {
\left[\begin{smallmatrix} 1&0\\ 0&1
 \end{smallmatrix}\right]}
\to\{\left[\begin{smallmatrix} 0&1\\
\lambda &0
 \end{smallmatrix}\right]\}_{\lambda
 \ne\pm 1}.\] These two arrows are also
correct since for each small
$\varepsilon\ne 0 $ we have
\begin{equation}\label{juw}
\begin{bmatrix}
0&-1+\varepsilon
\\
1 &1
 \end{bmatrix}\cong
\begin{bmatrix}
0&-1+\varepsilon
\\
1 &0
 \end{bmatrix}\cong
\begin{bmatrix}
0&1
\\
-1+\varepsilon &0
 \end{bmatrix}
\end{equation}
and
\[
\begin{bmatrix}
1&0
\\
\varepsilon&1
 \end{bmatrix}\cong
\begin{bmatrix}
0&1
\\
\lambda  &0
 \end{bmatrix}\quad\text{if }
2-\varepsilon^2=
\lambda+\lambda ^{-1} \text{ (see \eqref{htx}
and \eqref{ldt}).}
\]

(b)  The vertices of graph \eqref{g4}
represent the congruence bundles of
$3\times 3$ matrices, which are
described in Theorem \ref{theorm}.
Taking the structure of the graph $G_3$
in \eqref{g3} into account, it remains
to verify that each arrow with the
vertex $\{[\lambda
]^20\}_{\lambda\ne\pm 1}$ or $\{[\mu
]^21\}_{\mu \ne\pm 1}$ in \eqref{g4} is
correct and that we did not omit any
arrow with these vertices.

The arrows
\[
1^20\to \{[\lambda ]^20\}
_{\lambda\ne\pm 1},\quad
110\to \{[\lambda ]^20\}
_{\lambda\ne\pm 1},\quad
1^21\to \{[\mu ]^21\}
_{\mu\ne\pm 1}
\]
follow from the arrows
\[
{ \left[\begin{smallmatrix} 0&-1\\
1&1
 \end{smallmatrix}\right]}
\to\{\left[\begin{smallmatrix} 0&1\\
\lambda &0
 \end{smallmatrix}\right]\}_{\lambda
 \ne\pm 1},\quad
{
\left[\begin{smallmatrix} 1&0\\ 0&1
 \end{smallmatrix}\right]}
\to\{\left[\begin{smallmatrix} 0&1\\
\lambda &0
 \end{smallmatrix}\right]\}_{\lambda
 \ne\pm 1},\quad
 { \left[\begin{smallmatrix} 0&-1\\
1&1
 \end{smallmatrix}\right]}
\to\{\left[\begin{smallmatrix} 0&1\\
\lambda &0
 \end{smallmatrix}\right]\}_{\lambda
 \ne\pm 1}
\]
in \eqref{g2}.

The arrow $\{[\lambda
]^20\}_{\lambda\ne\pm 1}\to 0^3$
follows from the arrow $[\lambda
]^20\to 0^3$ in the graph $G_3$
constructed in \eqref{g3}.

Let us prove $1^3\to \{[\mu ]^21\}
_{\mu\ne\pm 1}$. Consider the following
perturbation of $1^3$:
\[
B_c:=\begin{bmatrix}
0 &0&1\\ 0&-1&-1\\c&1&0
 \end{bmatrix},\quad c\ne 1,
\]
in which $c$ is arbitrarily close to
$1$. The matrix $B_c$ is congruent to a
matrix from \eqref{g2v}. The symmetric
part of $B_c$ is nonsingular and its
skew-symmetric part is nonzero; the
only matrices in \eqref{g2v} with these
properties are $[\mu ]^21$ and $1^3$.
The eigenvalues of the cosquare $\Gamma
_3^{-T}\Gamma _3$ of $1^3$ are $1,1,1$
(see \eqref{kef}), and the eigenvalues
of the cosquare  $B_c^{-T}B_c$ are
$c,1,c^{-1}$. Thus, $ B_c\not\cong 1^3$
by Proposition \ref{leq}(b), and so
$B_c\cong [\mu ]^21$ for some $\mu $.

There are no arrows $\{[\lambda
]^20\}_{\lambda\ne\pm 1}\to [-1]^21$
and $\{[\lambda ]^20\}_{\lambda\ne\pm
1}\to 111$ since there are no directed
paths from $[\lambda ]^20$ to $[-1]^21$
and from $[\lambda ]^20$ to $111$ in
$G_3$  (see \eqref{g3}). There are no
arrows $[-1]^21\to \{[\lambda
]^20\}_{\lambda\ne\pm 1}$ and $111\to
\{[\lambda ]^20\}_{\lambda\ne\pm 1}$
since $[-1]^21$ and $111$ are
nonsingular, whereas all $[\lambda
]^20$ are
singular.\hfill\hbox{\qedsymbol}

\section{On the choice of definition
for congruence bundles}\label{kjr}

One could argue that Theorem
\ref{the1_b} about the closure graphs
$G_2^{\rm{B}}$ and $G_3^{\rm{B}}$ is of
doubtful value since these graphs may
change if another definition of
congruence bundles is used. Moreover,
we define congruence bundles of
matrices via bundles of pencils (see
Definition \ref{kuq1}) that may not
look natural since  similarity bundles
of matrices and bundles of matrix
pencils are defined via their canonical
forms.

However, Roiter's remarkable assertion
given in Lemma 4.1 partially justifies
Definition \ref{kuq1}. We give another
argument in support of Definition
\ref{kuq1}. We proceed from the
assumption that the partition of
$\mathbb C^{n\times n}$ into bundles
must satisfy at least the following two
conditions:
\begin{itemize}
  \item[(i)] each bundle is a union
      of congruence classes whose
      matrices have the same
      properties with respect to
      perturbations (which means
that the vertices representing
these congruence classes in the
closure graph $G_n$ are
indistinguishable in $G_n$), and

  \item[(ii)] each bundle consists
      of matrices that have the
      same number of indecomposable
      direct summands \eqref{cjy}
      in the congruence canonical
      form.
\end{itemize}
In Section \ref{kjy} we prove that the
partition of $\mathbb C^{n\times n}$
($n=2$ or $3$) into bundles in the
sense of Definition \ref{kuq1}
satisfies (i) and (ii); moreover, it is
the coarsest partition satisfying them.
Therefore, $G_2^{\rm{B}}$ and
$G_3^{\rm{B}}$ are correctly used as
the closure graphs. (In actuality, the
authors of this article first
constructed the closure graphs
$G_2^{\rm{B}}$ and $G_3^{\rm{B}}$ using
(i) and (ii), and then they defined
congruence bundles so as to get the
obtained $G_2^{\rm{B}}$ and
$G_3^{\rm{B}}$.)

In Section \ref{kjt} we show that the
congruence bundles do not satisfy (i)
if they are defined as the sets of
matrices whose congruence canonical
forms have the same type (i.e., the
canonical forms coincide up to their
sets of distinct parameters, in a
manner like \eqref{lic}).  In Section
\ref{kjtf} we find the isometry groups
of $2\times 2$ and $3\times 3$
congruence canonical matrices. In
Section \ref{kjt1} we discuss a reason
to not define a congruence bundle as a
set of matrices whose congruence
canonical forms have the same isometry
group.

\subsection{Arguments in support of
the definition of congruence
bundles}\label{kjy}

Let us show that our definition of
congruence bundles agrees with the
above condition (i) if $n=2$ or $3$.
This means that the partition of the
set of $2\times 2$ or $3\times 3$
matrices into congruence bundles is a
refinement of the partition into
classes of graph-equivalent matrices in
the sense of the following definition.

\begin{definition}\label{kda}
\begin{itemize}
  \item[\rm(a)] Two vertices $v$
      and $w$ of the closure
graph $G_n$ are
      \emph{graph-equivalent} if
      $v$
and $w$ are located at the same
horizontal level in $G_n$, and for
every vertex $x$
\[
v\to x\ \Longleftrightarrow\ w\to x,
      \qquad\qquad
v\leftarrow x\ \Longleftrightarrow\
w\leftarrow x.
\]

  \item[\rm(b)] Two $n\times n$
      matrices are
      \emph{graph-equivalent} if
      the vertices that correspond
      to their congruence classes
      are graph-equivalent.

\end{itemize}
\end{definition}
Thus, two vertices $v$ and $w$ of $G_n$
are graph-equivalent if and only if the
corresponding congruence classes have
the same dimension and the interchange
of $v$ and $w$ is an automorphism of
the directed graph $G_n$.

\begin{proposition}\label{juye}
\begin{itemize}
  \item[\rm(a)] The partition of
      $\C^{2\times 2}$ into
      congruence bundles is a
      refinement of the partition
into graph-equivalence classes.
The closure graph for
  graph-equivalence classes of\/
      $2\times 2$ matrices is given
      in Figure \ref{fig5}.
 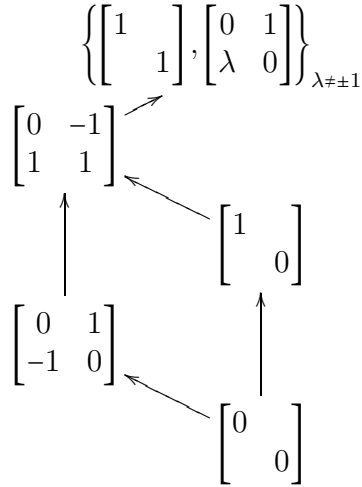
\begin{figure}[hbt]
\begin{equation*}\label{gke}
\begin{split}
\xymatrix@C=-20pt@R=1pt{
&{
\left\{\begin{bmatrix} 1&\\
&1 \end{bmatrix},\begin{bmatrix}
0&1\\ \lambda &0
 \end{bmatrix}\right\}_{\lambda \ne\pm 1}
 \hspace{30pt}}
\\ {\begin{bmatrix}
0&-1\\ 1&1
 \end{bmatrix}}\ar[ur]
 &
 &
    \\
&{\begin{bmatrix}
1&\\ &0
 \end{bmatrix}} \ar[lu]
    \\
{\begin{bmatrix} 0&1\\
-1&0
 \end{bmatrix}}\ar[uu] &
      \\
&{\begin{bmatrix} 0&\\
&0
 \end{bmatrix}}
\ar[ul]\ar[uu] }
\end{split}
\end{equation*}
\caption{\small The closure
  graph for graph-equivalence classes
      of\/ $2\times 2$ matrices (compare with the closure graph in Figure~\ref{fig3} for
congruence bundles of $2\times 2$ matrices).}
      \label{fig5}
\end{figure}
The vertex $\left\{\left[
\begin{smallmatrix} 1&\\
&1
 \end{smallmatrix}\right],
 \left[
\begin{smallmatrix} 0&1\\
\lambda &0
 \end{smallmatrix}\right]
 \right\}_{\lambda \ne \pm 1}$
represents the graph-equivalence
class consisting of all matrices
whose congruence canonical forms
are $\left[
\begin{smallmatrix} 1&\\
&1
 \end{smallmatrix}\right]$ or
$\left[\begin{smallmatrix} 0&1\\
\lambda &0
 \end{smallmatrix}\right]$ with $\lambda
 \ne \pm 1$. The other vertices are
canonical matrices under
congruence; their graph-equivalence
classes coincide with the
congruence classes.

  \item[\rm(b)] The partition of
      $\C^{3\times 3}$ into
      congruence bundles coincides
      with its partition into
graph-equivalence classes.
\hfill\hbox{\qedsymbol}

\end{itemize}
\end{proposition}

The graph-equivalence classes of
$2\times 2$ matrices coincide with the
congruence bundles except for one
graph-equivalence class that is the
union of two bundles: the bundle
containing
$\left[\begin{smallmatrix} 1&\\
&1
 \end{smallmatrix}\right]$ and the
bundle containing all
$\left[\begin{smallmatrix} 0&1\\
\lambda &0
 \end{smallmatrix}\right]$ with $\lambda
 \ne \pm 1$.
However, we can add a condition on
direct-sum-decompositions and obtain
partitions of the sets of $2\times 2$
and $3\times 3$ matrices that coincide
with their partitions into congruence
bundles:

\begin{definition}\label{kdah}
Two $n\times n$ matrices $A$ and $B$
are \emph{strictly graph-equivalent} if
\begin{itemize}
  \item[\rm(i)] $A$ and $B$ are
      graph-equivalent and
  \item[\rm(ii)] if $A$ and $B$ are
      congruent to
\begin{equation}\label{kke}
A_1\oplus\dots\oplus
A_k\qquad\text{and}\qquad
B_1\oplus\dots\oplus B_l,
\end{equation}
respectively, in
which each summand
      is a square matrix that is
      not congruent to a direct sum
      of square matrices of smaller
      sizes, then $k=l$.
\end{itemize}
\end{definition}

\begin{proposition}\label{jdw}
For $2\times 2$ and
      $3\times 3$ matrices, the
      congruence bundles coincide with
 the strict
      graph-equivalence
      classes.
      \hfill\hbox{\qedsymbol}
\end{proposition}

Proposition \ref{jdw} shows that the
definition of congruence bundles for
$2\times 2$ and $3\times 3$ matrices is
natural. It is an open problem whether
Proposition \ref{jdw} holds for
$n\times n$ matrices with $n>3$ even if
we extend (ii) in Definition~\ref{kdah}
by the condition: there is a suitable
renumbering of the summands in
\eqref{kke} such that each $A_i$ has
the same size as $B_i$ and they are
graph-equivalent in the set of matrices
of this size.

\subsection{A reason to not define a
congruence bundle as a set of matrices
having congruence canonical forms of
the same type} \label{kjt}

The bundles of matrices under
similarity and of matrix pencils under
strict equivalence are defined via the
Jordan and Kronecker canonical forms:
each bundle consists of all matrices or
pencils with the same canonical form
but with unspecified parameters; see
\eqref{lic}. In a similar way, one
could define bundles of matrices under
congruence as parameter-equivalence
classes in the sense of the following
definition.

\begin{definition}\label{mos}
Let $A$ be a square matrix. Write its
      congruence canonical form given
      in Proposition
      \ref{le1} as follows:
\begin{equation*}\label{lmso}
{\cal C}\oplus {\cal D}_1(\lambda_1)
\oplus\dots\oplus {\cal D}_t(\lambda_t),
\quad \text{$\lambda_i\ne \lambda _j,
 1/\lambda _j$ if }
i\ne j,\quad \lambda_1,\dots,\lambda_t
\in\mathbb C\setminus\{0\},
\end{equation*}
in which ${\cal C}$ is a direct sum of
matrices of the form $\Gamma _n$ and
$J_k(0)$, and each ${\cal
D}_i(\lambda_i)$ is a direct sum of
matrices of the form
$H_{2m}(\lambda_i)$ with $m$ such that
$\lambda_i\ne (-1)^{m+1}$. Then the
\emph{parameter-equivalence class} of
$A$ consists of all matrices whose
congruence canonical form is
\begin{equation*}\label{kfr}
{\cal C}\oplus {\cal D}_1(\mu_1)
\oplus\dots\oplus {\cal D}_t(\mu_t),
\quad \mu_i\ne\mu _j,1/ \mu _j\text{ if }
i\ne j,\quad \mu_1,\dots,\mu_t
\in\mathbb C\setminus\{0\},
\end{equation*}
with the same ${\cal C},{\cal D}_1,
\dots, {\cal D}_t$, in which $\mu_i\ne
(-1)^{m+1}$ for each direct summand
$H_{2m}(\mu_i)$ of ${\cal D}_i(\mu_i)$.
\end{definition}

We do not define congruence bundles as
parameter-equivalence classes (which
would be like bundles of matrices under
similarity and bundles of pencils)
since the matrices
\begin{equation}\label{hto}
\begin{bmatrix} 0&1\\
-1 &0
 \end{bmatrix}\qquad \text{and}\qquad
 \begin{bmatrix} 0&1\\
\lambda &0
 \end{bmatrix},\ \lambda
 \ne 0, \pm 1
\end{equation}
are parameter-equivalent, but they are
not graph-equivalent; see Remark
\ref{lrp}. Thus, the matrices
\eqref{hto} have distinct properties
with respect to perturbations and the
following proposition holds.

\begin{proposition}\label{jkf}
The partition of  $\C^{2\times 2}$ into
      parameter-equivalence classes is not a
      refinement of its partition
into graph-equivalence classes.
      \hfill\hbox{\qedsymbol}
\end{proposition}

Note that the matrices \eqref{hto}
belong to different bundles under
congruence, which is illustrated in
Figure~\ref{fig3}.

\subsection{The isometry
groups of $2\times 2$ and $3\times 3$
matrices}\label{kjtf}

The \emph{isometry group} of an
$n\times n$ matrix $A$ is the
multiplicative matrix group
\[
\st A:=\{S\in\mathbb C^{n\times n}\text{
is nonsingular}\,|\,
S^TAS=A\}
\]
(which is the isometry group of a
bilinear space $(\mathbb C^n,\mathcal
A)$, where $\cal A$ is the bilinear
form $(u,v)\mapsto u^TAv$ on $\mathbb
C^n$).

\DH okovi\'{c} \cite{dok} and Szechtman
\cite{sze} describe the structure of
$\st A$ for an arbitrary $A$. The group
$\st A$ has been much studied in the
case that $A$ is symmetric or skew
symmetric (see the references in
\cite{dok}); the theory of such groups
is close to the theory of classical Lie
groups. The following classical Lie
groups are isometry groups:
\begin{itemize}
  \item the general linear group
      $\gl_n(\C)=\st
      0_n=\{S\in\C^{n\times
      n}\,|\,\det S\ne 0\}$,
  \item the orthogonal group
      $\Or_n(\C)=\st
      I_n=\{S\in\C^{n\times
      n}\,|\,S^TS=I_n\}$,
  \item the symplectic group
      $\Sp_{2n}(\C)=\st \Omega
      _n=\{S\in\C^{n\times
      n}\,|\,S^T\Omega
      _nS=\Omega_n\}$,
\end{itemize}
in which
\[
\Omega_n:=\begin{bmatrix}
            0&I_n \\-I_n&0
          \end{bmatrix}.
\]
Note that
\[
\Sp_2(\C)=\SL_2(\C):=\{S\in\C^{2\times
      2}\,|\,\det S=1\}.\]

In the following theorem, we give  the
isometry groups of all $2\times 2$ and
$3\times 3$ congruence canonical
matrices (which represent the vertices
of $G_2$ and $G_3$ in \eqref{g1} and
\eqref{g3}). The isometry group of any
$2\times 2$ and $3\times 3$ matrix can
be obtained from the isometry group of
its congruence canonical matrix due to
the following assertion:
\[
\text{$B=R^TAR$ \ ($R$ is nonsingular)}
\quad\Longrightarrow\quad
\st B=R^{-1}(\st A)R.
\]
This assertion holds since for each
$S\in \st A$ we have $A=S^TAS$,
\[
R^TAR=R^TS^TR^{-T}\cdot R^TAR\cdot R^{-1}SR,
\]
$B=(R^{-1}SR)^T B (R^{-1}SR)$, and so
$R^{-1}SR\in \st B$.

\begin{theorem}\label{jlju}
\begin{itemize}
  \item[{\rm(i)}] The $2\times 2$
      canonical matrices under
      congruence have the following
      isometry
      groups:\\[2pt]
$\qquad\st \mat{0&\\
&0}= \gl_2(\mathbb C)$,
           \\[4pt]
$\qquad\st \mat{1&\\
&0} = \set{
\mat{\pm 1&0\\
*&c}}{c\ne 0} $,
\\[4pt]
$\qquad\st \mat{1&\\
&1} ={\rm O_2}(\mathbb C)$,
\\[4pt]
$\qquad\st \mat{0&1\\
-1&0} =\SL_2(\mathbb C)$,
\\[4pt]
$\qquad\st \mat{0&1\\
\lambda &0}=
 \left\{\left.\mat{c&0\\
0&1/c}\right| c\ne 0\right\}, \
\lambda \ne \pm 1$,
\\[4pt]
$\qquad\st \mat{0&-1\\
1&1} = \left\{
\pm\mat{1 &*\\
0&1}\right\},$ \\[2pt]
where the stars denote arbitrary
complex numbers.

 \item[{\rm(ii)}] The $3\times 3$
     canonical matrices under
     congruence have the following
     isometry
      groups:\\[2pt]
$\qquad\st \mat{0&&\\&0&\\&&0}=
\gl_3(\mathbb C)$,
\\[4pt]
$\qquad\st \mat{1&&\\&0&\\&&0}=
\set{ \mat{\!\pm 1\!
&\!\begin{smallmatrix}0&0
\end{smallmatrix}\!\\
\!\begin{smallmatrix}*\\ *
\end{smallmatrix}\!&\!B\!}}
{\det B\ne 0}$,
\\[4pt]
$\qquad\st \mat{1&&\\&1&\\&&0}=
\set{\mat{\!A
&\!\begin{smallmatrix}0\\0
\end{smallmatrix}\!\\
\!\begin{smallmatrix}*&*\end{smallmatrix}
&\!c\!}} {A^TA=I_2,\ c\ne 0}$,
\\[4pt]
$\qquad\st \mat{1&&\\&1&\\&&1}=
{\rm O_3}(\mathbb C)$,
\\[4pt]
$\qquad\st \mat{0&1&\\-1&0&\\&&0}=
\set{\mat{\!A
&\!\begin{smallmatrix}0\!\\0
\end{smallmatrix}\\\!
\begin{smallmatrix}*&*\end{smallmatrix}
&\!c\!}} {\det A=1,\ c\ne 0}$,
\\[4pt]
$\qquad\st \mat{0&1&\\-1&0&\\&&1}=
\set{\mat{\!A
&\!\begin{smallmatrix}0\!\\0
\end{smallmatrix}\\\!
\begin{smallmatrix}0&0\end{smallmatrix}
&\!\pm 1 \!}} {\det A=1}$,
\\[4pt]
$\qquad\st \mat{0&1&\\\lambda
&0&\\&&0}= \set{\mat{c&0&0\\
0&1/c&0\\ *&*&c'}} {c,c'\ne 0},\
\lambda \ne\pm 1$,
\\[4pt]
$\qquad\st \mat{0&1&\\\lambda
&0&\\&&1}=
\set{\mat{c&0&0\\ 0&1/c&0\\
0&0&\pm 1}} {c\ne 0},\ \lambda
\ne\pm 1$,
\\[4pt]
$\qquad\st \mat{0&-1&\\1 &1&\\&&0}=
\set{\mat{\delta &*&0\\ 0&\delta &0\\
*&*&c}} {\delta =\pm 1,\ c\ne 0}$,
\\[4pt]
$\qquad\st \mat{0&-1&\\1 &1&\\&&1}=
\set{\mat{\delta &*&0\\ 0&\delta &0\\
0&0&\pm 1}} {\delta=\pm 1}$,
\\[4pt]
$\qquad\st \mat{0&1&0\\0
&0&1\\0&0&0}=
\set{\mat{c&a&0\\ 0&1/c&0\\
0&-a&c}} {c\ne 0,\ a\text{ is
arbitrary}}$,
\\[4pt]
$\qquad\st \mat{0&0&1\\0
&-1&-1\\1&1&0}=
\set{\pm\mat{1&a&a^2/2\\
0&1&a\\
0&0&1}} {a\text{ is arbitrary}}$.
\end{itemize}
\end{theorem}

We omit the proof, which is a matter of
straightforward computation.

\subsection{A reason to not define
a congruence bundle as a set of
matrices whose congruence canonical
forms have the same isometry
group}\label{kjt1}

The \emph{stabilizer subgroup} of an
$n\times n$ matrix $A$ with respect to
similarity action is the multiplicative
matrix group
\[
G_A:=\{S\in\mathbb C^{n\times n}\text{
is nonsingular}\,|\,
S^{-1}AS=A\}.
\]
Two matrices $A$ and $B$ belong to the
same similarity bundle if and only if
the Jordan blocks of their Jordan
canonical forms $J(A)$ and $J(B)$ can
be arranged such that
$G_{J(A)}=G_{J(B)}$, which follows from
the description of the set of matrices
that commute with a Jordan matrix (see
\eqref{lic} and \cite[Chapter VIII,
\S\,2]{gan}).

Developing this approach, Patera,
Rousseau, and Schlomiuk \cite{pat} (see
also the references 1--3 in \cite{pat})
study partitions of complex and real
classical Lie algebras $\Lambda $ into
strata: they define a \emph{stratum} as
a set of matrices whose stabilizer
subgroups are conjugate (see \cite[p.
494]{pat}). They state that the
partition of $\Lambda =\mathbb
C^{n\times n}$ into similarity bundles
coincides with its partition into
strata and study the interrelations of
these partitions for other classical
Lie algebras (see \cite[p. 492]{pat}).

\begin{definition}\label{mowe}
We say that two $n\times n$ matrices
$A$ and $B$ are \emph{group-equivalent}
if the direct summands of their
congruence canonical forms $A\ca$ and
$B\ca$ can be arranged such that their
isometry groups coincide: $\st
A_{\text{\rm can}}= \st B_{\text{\rm
can}}$.
\end{definition}

\begin{proposition}\label{jle1}
\begin{itemize}
  \item[\rm(a)] For $2\times 2$ and
      $3\times 3$ matrices, the
      congruence bundles coincide
      with the group-equivalence
      classes.

  \item[\rm(b)] There are
      infinitely many
      group-equivalence classes
      of\/ $5\times 5$ matrices.
\end{itemize}
\end{proposition}

\begin{proof}
(a) The bundles of $2\times 2$ and
$3\times 3$ matrices were found in
Theorem \ref{the1_b}; they correspond
to the vertices of the closure graphs
$G_2^{\rm B}$ and $G_3^{\rm B}$. By
Theorem \ref{jlju}, all congruence
canonical matrices in a bundle have the
same isometry group, and any two
canonical matrices from different
bundles have unequal isometry groups.

(b) The matrices
\[
A_{\lambda }:=\left[\begin{array}{cc|ccc}
               0 & 1 &0 &0 &0 \\
               \lambda &0&0&0&0\\\hline
               0&0&0&1&0\\
               0&0&0&0&1\\
               0&0&0&0&0
             \end{array}
             \right]
\]
with distinct $\lambda\in\mathbb C$ are
not group-equivalent since for
\[
S_{\mu }:=\left[\begin{array}{cc|ccc}
               1& 0&0 &0 &0 \\
               0&1&0&1&0\\\hline
               -1&0&1&0&0\\
               0&0&0&1&0\\
               -\mu &0&0&0&1
             \end{array}
             \right]
\]
we have
\[
S_{\mu
}^TA_{\lambda }S_{\mu }= A_{\lambda }
\ (\text{i.e., }S_{\mu }\in \st(A_{\lambda }))
\quad\Longleftrightarrow\quad \lambda =\mu.
\tag*{\qedhere}
\]
\end{proof}

By Proposition \ref{jle1}(b), it would
be unwise to define congruence bundles
as group-equivalence classes.

\section*{Acknowledgement}

A.	Dmytryshyn and	B.	
K\r{a}gstr\"{o}m are supported by the
Swedish Research Council (VR) under
grant A0581501, and by eSSENCE, a
strategic collaborative e-Science
programme funded by the Swedish
Government via VR. V. Futorny is
supported by the CNPq (grant
301320/2013-6) and FAPESP (grant
2010/50347-9). The work of V.V.
Sergeichuk was done during a visit to
the University of S\~ao Paulo. He is
grateful to the University of S\~ao
Paulo for hospitality and to the FAPESP
for financial support (grant
2012/18139-2).

\end{document}